\documentclass{icmart}

\usepackage{amsmath,amssymb,amscd,mathrsfs,epic,wasysym,latexsym,tikz,mathrsfs,cite,hyperref}
\usepackage{pb-diagram}
\usepackage[matrix,arrow]{xy}

\contact[klesh@uoregon.edu]{Department of Mathematics, University of Oregon, Eugene, OR 97403, USA}

\newtheorem{Proposition}{Proposition}

\newtheorem{Theorem}[Proposition]{Theorem}
\newtheorem{Corollary}[Proposition]{Corollary}
\theoremstyle{definition}

\renewcommand{\theenumi}{\roman{enumi}}

\let\<\langle
\let\>\rangle

\newcommand\Comment[2][\relax]{\space\par\medskip\noindent%
   \fbox{\begin{minipage}{\textwidth}\textbf{Comment\ifx\relax#1\else---#1\fi}\newline%
        #2\end{minipage}}\medskip
}

\DeclareMathOperator\ledom{\trianglelefteq}

\def\bi{\text{\boldmath$i$}}
\def\bj{\text{\boldmath$j$}}

\def\b1{\text{\boldmath$1$}}

\def\pmod#1{\text{ }(\text{\rm mod } #1)\,}

\newcommand{\Hom}{\operatorname{Hom}}

\newcommand{\End}{\operatorname{End}}
\newcommand{\ind}{\operatorname{ind}}
\newcommand{\im}{\operatorname{im}}

\def\sgn{\mathtt{sgn}}
\newcommand{\res}{\operatorname{res}}
\newcommand{\soc}{\operatorname{soc}\,}
\newcommand{\head}{\operatorname{head}}

\newcommand{\St}{\operatorname{St}}
\newcommand{\infl}{\operatorname{infl}}

\newcommand{\Stab}{\operatorname{Stab}}

\newcommand{\Z}{\mathbb{Z}}

\def\eps{{\varepsilon}}
\def\phi{{\varphi}}

\newcommand{\ga}{\gamma}
\newcommand{\Ga}{\Gamma}
\newcommand{\la}{\lambda}
\newcommand{\La}{\Lambda}
\newcommand{\al}{\alpha}
\newcommand{\be}{\beta}

\newcommand{\Si}{\Sigma}
\newcommand{\si}{\sigma}

\newcommand{\de}{\delta}
\newcommand{\De}{\Delta}

\def\triv#1{\O_{#1}}

\def\Belt{\mathbf{B}}

\newcommand{\Mull}{{\bf M}}

\def\ImS{{\mathscr S}}

\newcommand{\Ind}{{\mathrm {Ind}}}

\newcommand{\tr}{{\mathrm {tr}}}

\newcommand{\pr}{{\mathrm {pr}}}

\newcommand{\Res}{{\mathrm {Res}}}
\newcommand{\Ann}{{\mathrm {Ann}}}
\newcommand{\C}{{\mathbb C}}
\newcommand{\Q}{{\mathbb Q}}

\newcommand{\A}{{\mathscr A}}
\newcommand{\D}{{\mathscr D}}

\renewcommand{\mod}{\bmod \,}

\newcommand{\cont}{\operatorname{cont}}

\def\Par{{\mathscr P}}
\def\ula{{\underline{\lambda}}}
\def\umu{{\underline{\mu}}}
\def\unu{{\underline{\nu}}}

\def\JTD{{\mathbf D}}

\def\g{\mathfrak{g}}

\def\GG{\text{\boldmath$g$}}

\def\T{{\mathtt T}}

\def\Stab{{\mathtt S}}
\def\G{{\mathtt G}}

\def\spa{\operatorname{span}}
\def\height{{\operatorname{ht}}}

\def\wt{{\operatorname{wt}}}

\def\RPar{{\mathscr{RP}}}
\def\surj{{\twoheadrightarrow}}
\def\Gar{{\operatorname{Gar}}}
\def\op{{\mathrm{op}}}
\def\re{{\mathrm{re}}}
\def\im{{\mathrm{im}\,}}

\def\mod#1{#1\!\operatorname{-mod}}
\def\proj#1{#1\!\operatorname{-proj}}

\def\HCI{I}

\def\iso{\stackrel{\sim}{\longrightarrow}}

\def\Mde{M}

\def\Zde{Z}
\def\Lade{\La}

\def\Zdot{Z}
\def\Lde{L}
\def\Dede{\De}
\def\nade{\nabla}

\def\triv{{\tt triv}}

\def\Stand{\Delta}

\def\CH{{\operatorname{ch}_q\,}}
\def\DIM{{\operatorname{dim}_q\,}}

\def\words{{\langle I\rangle}}%{{\mathbf W}}
\def\shift{{\tt sh}}

\def\Seq{{\tt Se}}
\def\Car{{\tt C}}
\def\cc{{\tt c}}

\newcommand{\Laurent}{\mathcal L}

\newcommand\SetBox[2][35mm]{\Big\{\vcenter{\hsize#1\centering#2}\Big\}}

{\catcode`\|=\active
  \gdef\set#1{\mathinner{\lbrace\,{\mathcode`\|"8000%
  \let|\midvert #1}\,\rbrace}}
}
\def\midvert{\egroup\mid\bgroup}

\colorlet{darkgreen}{green!50!black}
\tikzset{dots/.style={very thick,loosely dotted},
         greendot/.style={fill,circle,color=darkgreen,inner sep=1.5pt,outer sep=0}
}
\def\greendot(#1,#2){\node[greendot] at(#1,#2){}}

\newenvironment{braid}{% sets defaults for the braid diagrams
  \begin{tikzpicture}[baseline=6mm,blue,line width=1pt, scale=0.4,
                      draw/.append style={rounded corners},
                      every node/.append style={font=\fontsize{5}{5}\selectfont}]%
  }{\end{tikzpicture}
}

\def\Grid(#1,#2){%  draws a coordinate grid inside a braid diagram
  \draw[very thin,gray,step=2mm] (0,0)grid(#1,#2);
  \draw[very thin,darkgreen,step=10mm] (0,0)grid(#1,#2);
}

\newcommand\Tableau[2][\relax]{
  \begin{tikzpicture}[scale=0.5,draw/.append style={thick,black}]
    \ifx\relax#1\relax%
    \else % shade the boxes in #1
      \foreach\box in {#1} { \filldraw[blue!30]\box+(-.5,-.5)rectangle++(.5,.5); }
    \fi
    \newcount\row\newcount\col
    \row=0
    \foreach \Row in {#2} {
       \col=1
       \foreach\k in \Row {
          \draw(\the\col,\the\row)+(-.5,-.5)rectangle++(.5,.5);
          \draw(\the\col,\the\row)node{\k};
          \global\advance\col by 1
       }
       \global\advance\row by -1
    }
  \end{tikzpicture}
}

\newcommand\YoungDiagram[2][\relax]{
  \begin{tikzpicture}[scale=0.5,draw/.append style={thick,black}]
    \ifx\relax#1\relax%
    \else % shade the boxes in #1
    \foreach\box in {#1} {
      \filldraw[blue!30]\box rectangle ++(1,1);
    }
    \fi
    \newcount\row
    \row=0
    \foreach \col in {#2} {
       \draw(1,\the\row)grid ++(\col,1);
       \global\advance\row by -1
    }
  \end{tikzpicture}
}

\title[Representations of Symmetric Groups]{Modular Representation Theory of Symmetric Groups}

\author[Alexander Kleshchev]
{Alexander Kleshchev\thanks{Research 
supported by the NSF grant no. DMS-1161094 and the Humboldt Foundation.}}

\begin{document}

\begin{abstract}
We review some recent advances in modular representation theory of symmetric groups and related Hecke algebras. We discuss  connections with Khovanov-Lauda-Rouquier algebras and  gradings on the blocks of the group algebras $F\Si_n$, which these connections reveal; graded categorification and connections with quantum groups and crystal bases; modular branching rules and the Mullineaux map; graded cellular structure and graded Specht modules; cuspidal systems for affine KLR algebras and imaginary Schur-Weyl duality, which connects representation theory of these algebras to the usual Schur algebras of smaller rank. 
\end{abstract}

\begin{classification}
Primary 20C30; Secondary 20C08, 17B37.
\end{classification}

\begin{keywords}
Symmetric group, modular representation, Hecke algebra, quantum group. 
\end{keywords}

\maketitle

\section{Introduction}
The classical problem of understanding representation theory of symmetric groups is especially difficult in positive characteristic. For example, there is no effective algorithm for computing the dimensions of irreducible modules in that case. 
In this expository paper, we review some results on {\em modular}\,  representation theory of symmetric groups. Let $F$ be a filed of characteristic $p> 0$ and $\Si_n$ be the symmetric group on $n$ letters. 

We discuss connections with Khovanov-Lauda-Rouquier algebras and gradings on the blocks of $F\Si_n$ which these connections  reveal, graded categorification and connections with quantum groups and crystal bases, modular branching rules and the Mullineaux map, graded cellular structure and graded Specht modules, cuspidal systems for affine KLR algebras and imaginary Schur-Weyl duality which connects representation theory of these algebras to the usual Schur algebras of smaller rank. 

\section{Graded isomorphism theorem}\label{SGI}

\subsection{Basic Notation}
%Throughout the paper we assume that the characteristic $p$ of the ground field $F$ is positive. 
Let $I:=\Z/p\Z$ be identified with the simple subfield of $F$. We associate to $p$ the affine {\em Cartan matrix}\, $\Car=(\cc_{ij})_{i,j\in I}$ of type ${\tt A_{p-1}^{(1)}}$ so that $\cc_{ii}=2, \cc_{ij}=-1$ if $|i-j|=1$ and $p\neq 2$, $\cc_{0,1}=\cc_{1,0}=-2$ if $p=2$ and $\cc_{ij}=0$ otherwise. 

\iffalse{
Let $\Gamma$ be the cyclic graph with vertices $i\in I$ connected if and only if $|i-j|=1$. 
For example, if $p=5$:  

$$\begin{tikzpicture}%[scale=4.5, line join=bevel]
\node at (-1.7,1) {$\Ga=$};
\node at (0,2) {$0$};
\node at (1,1.3) {$1$}; 
\node at (0.5,0.3) {$2$}; 
\node at (-0.5,0.3) {$3$}; 
\node at (-1,1.3) {$4$}; 
\draw [-] (0.85,1.35) -- (0.15,1.85);
\draw [-]  (0.6,0.5) -- (0.87,1.1);
\draw [-]  (-0.3,0.25) -- (0.3,0.25);
\draw [-]  (-0.87,1.1) -- (-0.6,0.5);
\draw [-] (-0.15,1.85) -- (-0.85,1.35);
\end{tikzpicture}
$$
}\fi

The group algebra $F\Si_n$ will be denoted by $H_n$. 
We have the simple transpositions $s_r:=(r,r+1)\in\Si_n$ for $r=1,\dots n-1$, and the {\em Murphy elements}: 
$$
x_t:=(1,t)+(2,t)+\dots+(t-1,t)\in H_n\qquad(1\leq t\leq n). 
$$
The Murphy elements commute. 

Let $V$ be a finite dimensional $H_n$-module. 
The eigenvalues of the Murphy elements in $V$ belong to $I$, see for example \cite[Lemma 7.1.2]{Kbook}. Given a {\em word} $\bi=(i_1,\dots,i_n)\in I^n$, we define the
{\em $\bi$-word space} of $V$ as follows:
$$
V_\bi=\{v\in V\mid (x_r-i_r)^Nv=0\ \text{for $N\gg0$ and $r=1,\dots,n$}\}.
$$
We have a {\em word space decomposition}\, 
$
V=\bigoplus_{\bi\in I^n} V_\bi.
$
Using the word space decomposition of the left regular $H_n$-module, we get a system of orthogonal idempotents 
$ %\begin{equation}\label{EIdempotents}
\{1_\bi\mid \bi\in I^n\}
$ %\end{equation}
in the group algebra $H_n$ (some of which are zero) such that
$
\sum_{\bi\in I^n}1_\bi=1,
$
and 
$
1_\bi V=V_\bi 
$
for all $\bi\in I^n$ and all finite dimensional $H_n$-modules~$V$.

The symmetric group $\Si_n$ acts on $I^n$ by place permutations, and let $I^n/\Si_n$ be the set of orbits. 
Fix an orbit $\al\in I^n/\Si_n$. Define 
$%\begin{equation}\label{bal}
1_{\alpha} := \sum_{\bi \in \al} 1_\bi \in H_n.
$ %\end{equation}
It is easy to check that $1_\al$ is a central idempotent in $H_n$. In fact, by \cite{LM} or \cite[Theorem~1]{cyclo},
$1_{\alpha}$ is either zero or a primitive central idempotent
in $H_n$. 
Hence the algebra 
\begin{equation*}\label{fe}
H_\alpha := 1_{\alpha} H_n
\end{equation*}
is either zero or  a single {\em block}
of the group algebra $H_n$. If $x\in H_n$ we denote $1_\al x\in H_\al$ again by $x$. 

\subsection{Graded presentation}
Define special elements of $H_\al$ as follows: 
\begin{align}\label{QEPolKL}
y_t&:=\sum_{\bi\in \al}(x_t-i_r) 1_\bi &(1\leq t\leq n),
\\
\label{QQCoxKL}
\psi_r&:=
\sum_{\bi\in \al}(s_r+P_r(\bi))Q_r(\bi)^{-1}1_\bi &(1\leq r< n), \end{align}
where $P_r(\bi)$ and $Q_r(\bi)^{-1}$ are certain explicit polynomials in $F[y_r,y_{r+1}]$ defined
in~\cite{BK}. This gives us the following elements of $H_\al$:
\begin{equation}\label{EKLElts}
\{1_\bi\mid \bi\in \al\}\cup\{y_1,\dots,y_n\}\cup\{\psi_1,\dots,\psi_{n-1}\}.
\end{equation}

Finally, choose signs $\eps_{ij}\in{\pm1}$ for all $i,j \in I$ with $|i-j|=1$  so that $\eps_{ij}\eps_{ji} = -1$, and define the polynomials in $F[u,v]$:
\begin{equation}\label{EArun}
Q_{ij}(u,v):=
\left\{
\begin{array}{ll}
0 &\hbox{if $i=j$;}\\
1 &\hbox{if $\cc_{ij}=0$;}\\
%-u^{-\cc_{ij}}+v^{-\cc_{ji}} &\hbox{if $\cc_{ij}<0$ and $i> j$;}\\
\eps_{ij}(u^{-\cc_{ij}}-v^{-\cc_{ji}}) &\hbox{if $\cc_{ij}<0$.}
\end{array}
\right.
\end{equation}

The following result was first proved in \cite[Theorem 1.1]{BK}, see also \cite{R}: 

\begin{Theorem} \label{TBK}  %{\bf ()}
The algebra $H_\al$ is generated by the elements (\ref{EKLElts}) subject {\em only} to the following relations (for all admissible $\bi=(i_1,\dots),\bj, r,t$):
\begin{align}\label{R1}
1_\bi  1_\bj &= \de_{\bi,\bj} 1_\bi ,
\\
\quad{\textstyle\sum_{\bi \in \alpha}} 1_\bi  &= 1;\\
%\begin{equation}\label{R2PsiY}
y_r 1_\bi  &= 1_\bi  y_r;
\\
\quad y_r y_t &= y_t y_r;
\\
%\end{equation}
%\begin{equation}
\psi_r 1_\bi  &= 1_{s_r\bi} \psi_r;\label{R2PsiE}
\\
%\end{equation}
%\begin{equation}
(y_t\psi_r-\psi_r y_{s_r(t)})1_\bi  
&= \de_{i_r,i_{r+1}}(\de_{t,r+1}-\de_{t,r})1_\bi;
\label{R6}
\\
%\end{equation}
%\begin{equation}
\psi_r^21_\bi  &= Q_{i_r,i_{r+1}}(y_r,y_{r+1})1_\bi; 
 \label{R4}
\\
%\end{equation}
%\begin{equation} 
\psi_r \psi_t &= \psi_t \psi_r\qquad (|r-t|>1);\label{R3Psi}
\\
%\end{equation}
%\begin{equation}
%\begin{split}
(\psi_{r+1}\psi_{r} \psi_{r+1}-\psi_{r} \psi_{r+1} \psi_{r}) 1_\bi  
&=
\de_{i_r,i_{r+2}}\textstyle \frac{Q_{i_r,i_{r+1}}(y_{r+2},y_{r+1})-Q_{i_r,i_{r+1}}(y_r,y_{r+1})}{y_{r+2}-y_r}1_\bi;
%\end{split}
\label{R7}
\\
%\end{equation}
%\end{align}
%\begin{equation}
\label{ECyclotomic}
y_1^{\de_{i_1,0}}1_\bi&=0.
\end{align}
\end{Theorem}

We note that the ratio in (\ref{R7}) is always a polynomial in $y$'s. 
Theorem~\ref{TBK} is saying that blocks of groups algebras of  symmetric groups are isomorphic to certain {\em cyclotomic Khovanov-Lauda-Rouquier (KLR) algebras} \cite{KL1,KL2,R,R2}. We return to this in Section~\ref{SKLR}. 
The presentation of $H_\al$ given in Theorem~\ref{TBK} allows us to define a {\em grading} on $H_\al$ by setting:
$$
\deg(1_\bi ):=0,\quad \deg(y_r1_\bi ):=2,\quad\deg(\psi_r 1_\bi ):=-\cc_{{i_r},i_{r+1}}.
$$

From now on, we will always consider {\em graded} $H_\al$-modules, unless otherwise stated. The irreducible ungraded $H_\al$-modules are gradable in a unique way up to isomorphism and degree shift, so considering graded modules does divert our attention from the main goal of understanding irreducible $H_n$-modules. 

\subsection{Basics of graded algebra}
For any graded algebra $H$, we denote by $\mod{H}$ the abelian category of all finitely generated graded $H$-modules, with morphisms $\hom_H(\cdot,\cdot)$ being {\em degree-preserving} module homomorphisms. Denote by 
%$\mod{H}$ the full subcategory of finitely generated graded $H$-modules, and by 
$\proj{H}$ the full subcategory of finitely generated projective graded $H$-modules. 
%{\em finite dimensional}\, graded $H$-modules, with morphisms $\hom_H(\cdot,\cdot)$ being {\em degree-preserving} module homomorphisms. 
%We usually work in $\mod{H}$. 
Set
$$
\Laurent:=\Z[q,q^{-1}]. 
$$

The Grothendieck group $[\mod{H}]$ is a $\Laurent$-module via 
%\begin{equation}\label{EAMod}
$
q^m[M]:=[q^m M],
$ 
%\end{equation}
where $q^m M$ denotes the module obtained by 
shifting the grading in $M$ up by $m$, i.e. 
$
%\begin{equation}\label{obvious}
(q^mM)_n:=M_{n-m}.
%\end{equation}
$ 
%Given $f = \sum_{n \in \Z} f_n q^n \in \Z[[q,q^{-1}]]$ and $M\in\Mod{H}$,we write $$f \cdot M:=\bigoplus_{n \in \Z} M \langle n \rangle^{\oplus f_n}.$$
For $n \in \Z$, let
$
\Hom_H(M, N)_n := \hom_H(q^nM, N)
$, and 
set
$$
\Hom_H(M,N) := \bigoplus_{n \in \Z} \Hom_H(M,N)_n. 
$$ 
If $M$ is finitely generated, forgetting all gradings, $\Hom_H(M,N)$ is the usual Hom. 

For graded $H$-modules $M$ and $N$ we write $M\cong N$ to mean that $M$ and $N$ are isomorphic as graded modules and $M\simeq N$ to mean that $M\cong q^d N$ for some $d\in\Z$. 
%they are isomorphic as $H$-modules after we forget the gradings. 
For a finite dimensional 
graded vector space $V=\oplus_{d\in \Z} V_n$, its {\em graded dimension} is 
$$\DIM \, V:=\sum_{d \in \Z}  (\dim V_d)q^d\in\Laurent.$$ 
Given $M, L \in \mod{H}$ with $L$ irreducible, 
we write $[M:L]_q$ for the corresponding {\em  graded composition multiplicity},
i.e. 
%\begin{equation}\label{EGCM}
$
[M:L]_q := \sum_{n \in \Z} a_d q^d,
$ 
%\end{equation}
where $a_d$ is the multiplicity
of $q^dL$ in a graded composition series of $M$. 

Define the {\em formal character} of $M\in\mod{H_\al}$ as the formal sum 
$$%\begin{equation*}%\label{grch}
\CH M:=\sum_{\bi\in\al}(\DIM M_\bi) \bi. 
$$ %\end{equation*}
Consider the antiautomorphism 
\begin{equation}\label{EPrime}
H_\al\to H_\al,\ h\mapsto h',
\end{equation} 
which is identity on the generators (\ref{EKLElts}). If $M=\bigoplus_{d\in\Z}M_d\in\mod{H_\al}$, then the {\em graded dual}
$M^\circledast$ is the graded $H_\al$-module such that $(M^\circledast)_d:=M_{-d}^*$, for all
$d\in\Z$, and the $H_\al$-action is given by $(xf)(m)=f(x'm)$, for all $f\in M^\circledast, m\in M, x\in
H_\al$.  

For every irreducible $H_\al$-module $L$, there is a unique choice of the grading shift so that $L^\circledast \cong L$ \cite[Section~3.2]{KL1}. We always choose the shifts for irreducible $H_\al$-modules in this way. 

\subsection{The KLR algebra $A_\al$}\label{SSKLRDef}
We denote by $A_\al$ the algebra given by the generators (\ref{EKLElts}) and the relations (\ref{R1})--(\ref{R7}). This is the {\em KLR algebra} corresponding to $\Car$. It is an infinite dimensional graded algebra with the natural surjection 
\begin{equation}\label{ESurj}
A_\al\surj H_\al.
\end{equation}
All $H_\al$-modules will be considered as $A_\al$-modules via the functor of inflation $\infl:\mod{H_\al}\to \mod{A_\al}$. This functor has a left adjoint $\pr:\mod{A_\al}\to \mod{H_\al},\ M\mapsto H_\al\otimes_{A_\al} M$. The definition of the graded duality `$\circledast$'  and the formal character $\CH$ for $H_\al$-modules extends to $A_\al$-modules. 

In the important special case $\al=n\al_i$, the algebra $A_{n\al_i}$ is the usual affine nilHecke algebra. 
It has a representation on the polynomial space  $P_n=F[x_1,\dots,x_n]$
with each $y_t$ acting as multiplication by $x_t$ and each $\psi_r$ acting 
as the divided difference operator
$
f \mapsto ({^{s_r}} f - f)/(y_{r}-y_{r+1}).
$ 
The module $P_n$ is graded so that $\deg(x_r)=2$. 
We shift the degrees %of $P_n$ 
and define 
\begin{equation}\label{EPi^n}
P(i^{(n)}):=q^{-n(n-1)/2}P_n.
\end{equation}

%\vspace{4mm}\subsection{Specht modules James classification}

\section{Branching and categorification}
\subsection{More notation}%First, a little more notation. 
Following \cite[\S 1.1]{Kac}, we have a %let $(\h,\Pi,\Pi^\vee)$ be a 
realization of the Cartan matrix $\Car$. In particular we have the simple roots $\{\al_i\mid i\in I\}$, the fundamental dominant weights $\{\La_i\mid i\in I\}$, and the  form $(\cdot,\cdot)$ such that
$$
(\al_i,\al_j)=\cc_{i,j}\quad \text{and} \quad (\La_i,\al_j)=\de_{i,j}\qquad(i,j\in I). 
$$
Denote 
$Q_+ := \bigoplus_{i \in I} \Z_{\geq 0} \al_i$ and $P:=\bigoplus_{i\in I}\Z \La_i$. For $\alpha=\sum_{i\in I}n_i\al_i \in Q_+$, we write $\height(\alpha):=\sum_{i\in I} n_i$. %for the sum of its coefficients when expanded in terms of the $\al_i$'s. 
We have a bijection 
$$I^n/\Si_n\iso \{\al\in Q_+\mid \height(\al)=n\},\ \Si_n\cdot (i_1,\dots,i_n) \mapsto \al_{i_1}+\dots+\al_{i_n},$$ 
and from now on we {\em identify the two sets}. So we have the algebras $H_\al$ for $\al\in Q_+$.

%, the fundamental dominant weights $\{\La_i\mid i\in I\}$, and the normalized invariant form $(\cdot,\cdot)$ such that
%$$(\al_i,\al_j)=\cc_{i,j}, \quad (\La_i,\al_j)=\de_{i,j}\qquad(i,j\in I).$$

%Let $\g'=\g(\Car')$ be the finite dimensional simple Lie algebra  whose Cartan matrix $\Car'$ corresponds to the subset of vertices $I':=I\setminus\{0\}$. The affine Lie algebra $\g=\g(\Car)$ is then obtained from $\g'$ by a procedure described in \cite[Section 7]{Kac}. We denote by $W$ (resp. $W'$) the corresponding {\em affine Weyl group} (resp. {\em finite Weyl group}). It is a Coxeter group with standard generators $\{r_i\mid i\in I\}$ (resp. $\{r_i\mid i\in I'\}$), see \cite[Proposition~3.13]{Kac}. 

\subsection{Induction and restriction functors}
We want to study the induction and restriction functors between $H_n$-modules and $H_{n-1}$-modules. In particular, we are interested to know as much as possible about restrictions of irreducible $H_n$-modules to $H_{n-1}$, since this could help us to understand  irreducible $H_n$-modules by induction. 

It makes sense to refine induction and restriction to blocks. 
For any $\al\in Q_+$ of height $n$ and $i \in I$, there is an obvious graded algebra homomorphism
$H_\alpha \rightarrow H_{\alpha+\alpha_i}$. It 
maps the identity element of $H_\alpha$ to the idempotent
$$1_{\alpha,\alpha_i}= \sum_{\bi \in \alpha+\alpha_i,\ 
i_{n+1}=i} 1_\bi \in H_{\alpha+\alpha_i}.$$ Now, define the functors 
\begin{align*}
e_{i}&:=1_{\alpha,\alpha_i} H_{\alpha+\alpha_i} \otimes_{H_{\alpha+\alpha_i}} -
:\mod{H_{\alpha+\alpha_i}} \rightarrow \mod{H_{\alpha}},
\\
f_{i}&:=H_{\alpha+\alpha_i} 1_{\alpha,\alpha_i} \otimes_{H_{\alpha}} -
:\mod{H_\alpha} \rightarrow \mod{H_{\alpha+\alpha_i}}.
\end{align*}
\iffalse{
Let $L$ be an irreducible $H_n$-module. Then $L$ is in a single block of $H_n$, i.e. there is a unique $\al\in Q_+$ with $1_\al L\neq 0$, and $L$ can be considered as an irreducible $H_\al$-module. More generally, given any $M\in\mod{H_\al}$, considering it first as an $H_n$-module and then restricting to $H_{n-1}$, we get an $H_{n-1}$-module which does not necessarily belong to a single block. However, it is not difficult to see that $1_\be\res^{H_n}_{H_{n-1}}M\neq 0$ only if $\be$ is of the form $\be=\al-\al_i$ for some $i\in I$. 
Similarly $1_\be\ind_{H_n}^{H_{n+1}}M\neq 0$ only if $\be$ is of the form $\be=\al+\al_i$ for some $i\in I$. 
Define
$$
e_i M:= 1_{\al-\al_i}\res^{H_n}_{H_{n-1}}M, \quad
f_i M:= 1_{\al+\al_i}\ind_{H_n}^{H_{n+1}}M
\qquad(i\in I).
$$
This defines functors
$$
e_i : \mod{H_\al}\to \mod{H_{\al-\al_i}},\quad
f_i : \mod{H_\al}\to \mod{H_{\al+\al_i}}\qquad(\al\in Q_+,\ i\in I). 
$$
The functors extend by additivity to $\mod{H_n}$, so that for any $M\in\mod{H_n}$ we have
 $\res^{H_n}_{H_{n-1}}M \cong \bigoplus_{i\in I} e_i M$ and $\ind_{H_n}^{H_{n+1}}M \cong \bigoplus_{i\in I} f_i M$. 
}\fi 
For $M\in\mod{H_\al}$, define 
$$
\eps_i(M):=\max\{k\mid e_i^k\neq 0\},\quad \phi_i(M):=\max\{k\mid f_i^k\neq 0\}.
$$

\subsection{First branching rules}
Let %$\al\in Q_+$ and 
$L$ be an irreducible $H_\al$-module. 
It has been first proved in \cite{KBrIII} (in the ungraded setting) that $e_i L$ is either zero or it has a simple socle and head isomorphic to each other, and similarly for $f_iL$, see also \cite{GrVaz} for a more conceptual proof and a generalization. In particular, $e_iL$ and $f_iL$ are either zero or {\em indecomposable}, which is far from obvious. Let 
$$
\tilde e_i L:=\soc e_i L,\quad \tilde f_i L:=\soc f_i L. 
$$
This defines maps $$\tilde e_i, \tilde f_i:B\to B\sqcup\{0\},$$ where $B$ is the set of irreducible $H_\al$-modules for all $\al\in Q_+$ up to isomorphism and degree shift. Recall that all algebras and modules are graded, moreover the irreducible modules are graded {\em canonically} so that they are gradedly self-dual. This applies in particular to the irreducible modules $\tilde e_i L, \tilde f_i L$. 

For $n\in\Z$, denote the corresponding {\em quantum integer} 
$$[n]_q:=(q^n-q^{-n})/(q-q^{-1}).$$ Then one can refine the results on the socle as follows \cite[Theorem 4.12]{BKllt}: 

\begin{Theorem} \label{TBr}%{\rm \cite{}}%{\bf ()}
Let $\al\in Q_+$, $i\in I$ and $L$ be an irreducible $H_\al$-module. Then:
\begin{enumerate}
\item[{\rm (i)}] $(e_i L)^\circledast\cong e_i L$ and $(f_i L)^\circledast\cong f_iL$.
\item[{\rm (ii)}] $e_i L$ and $f_iL $ are indecomposable or zero. Moreover:  
\begin{align*}
\label{}
&\soc e_iL\cong q^{\eps_i(L)-1}\tilde e_iL, \quad &\head e_iL\cong q^{1-\eps_i(L)}\tilde e_iL,
\\
&\soc f_iL\cong q^{\phi_i(L)-1}\tilde f_iL,\quad &\head f_iL\cong q^{1-\phi_i(L)}\tilde f_iL.
\end{align*}
\item[{\rm (iii)}] $[e_i L: \tilde e_iL]_q=[\eps_i(\la)]_q$ and $[f_i L: \tilde f_iL]_q=[\phi_i(\la)]_q$. 
\item[{\rm (iv)}] $\eps_i(\tilde e_iL)=\eps_i(L)-1$ and $\eps_i(N)<\eps_i(L)-1$ for any other composition factor $N$ of $e_iL$; $\phi_i(\tilde f_iL)=\phi_i(L)-1$ and $\phi_i(K)<\phi_i(L)-1$ for any other composition factor $K$ of $f_iL$.  
\item[{\rm (v)}] $\End_{H_{\al-\al_i}}(e_iL)\cong F[x]/(x^{\eps_i(L)})$, the truncated polynomial algebra with the variable $x$ of degree $2$, and   $\End_{H_{\al+\al_i}}(f_iL)\cong F[x]/(x^{\phi_i(L)})$. 
\item[{\rm (vi)}] $e_iL$ is irreducible if and only if $\eps_i(L)=1$. \end{enumerate} 
\end{Theorem}

Ungraded versions of these results were first obtained in \cite{KBrII}, \cite{KBrIII}, \cite{KDec}, \cite{KBrIV}. For more branching rules see \S\ref{SSMoreBR}. 

\subsection{Crystal operators}
Let us return to the set $B$ and the operations $\tilde e_i,\tilde f_i:B\to B\sqcup\{0\}$. Every elements of $[L]\in B$ is an isomorphism class of an irreducible $H_\al$-module $L$ for some $\al\in Q_+$. This allows us to define a function $$\wt:B\to P,\ [L]\mapsto \La_0-\al.$$ Moreover, for all $i\in I$, we have functions $\eps_i,\phi_i:B\to \Z_{\geq 0}$.

Let $\g$ be the affine Kac-Moody Lie algebra corresponding to the Cartan matrix~$\Car$, i.e. $\g=\widehat {\mathfrak sl}_p(\C)$, see \cite[Section 7]{Kac}.

\begin{Theorem} \label{TCryst}%{\rm \cite{}}%{\bf ()}
The tuple $(B,\eps_i,\phi_i,\tilde e_i, \tilde f_i,\wt)$ is the Kashiwara's crystal associated to the irreducible $\g$-module $V(\La_0)$ with  highest weight $\La_0$. 
\end{Theorem}

This theorem has been first proved by Lascoux, Leclerc and Thibon \cite{LLT} by comparing the branching rules from \cite{KBrII} and \cite{KBrIII} with the explicit combinatorial description of the crystal obtained in \cite{MisMiw}. A more conceptual proof was found in \cite{Gr}, see also \cite{Kbook}. The remaining results of this section can be considered as steps towards an `explanation' of the theorem, the `real explanation' coming perhaps from \cite{VV} and \cite{R2}.

\subsection{Divided powers}
We define divided power analogues of the functors $e_i,f_i$. In order to do this, we exploit the algebras $A_\al$. There is an obvious embedding $A_{\al,\be}:=A_{\alpha}\otimes A_{\beta}\to A_{\alpha+\beta}$
mapping $1\otimes 1\mapsto 1_{\alpha,\beta} := \sum_{\bi\in \alpha,\, \bj \in \beta} 1_{\bi\bj}.
$ 
Consider the functors
\begin{align}
\Ind_{\alpha,\beta}^{\alpha+\beta} &:= A_{\alpha+\beta} 1_{\alpha,\beta}
\otimes_{A_{\alpha,\beta}} -:\mod{A_{\alpha,\beta}} \rightarrow \mod{A_{\alpha+\beta}},
\label{EIndFunctor}
\\
\label{EResFunctor}
\Res_{\alpha,\beta}^{\alpha+\beta} &:= 1_{\alpha,\beta} A_{\alpha+\beta}
\otimes_{A_{\alpha+\beta}} -:\mod{A_{\alpha+\beta}}\rightarrow \mod{A_{\alpha,\beta}}.
\end{align}
Let $i \in I$ and $n \geq 1$. Recalling the $A_{n\al_i}$-module $P(i^{(n)})$ from (\ref{EPi^n}), set
\begin{align*}%\label{div1}
\theta_{i}^{(n)}&:= 
\Ind_{\alpha,n \alpha_i}^{\alpha+n\alpha_i} (- \boxtimes P(i^{(n)})):\mod{A_\alpha} \rightarrow \mod{A_{\alpha+n\alpha_i}},\\
(\theta_{i}^*)^{(n)}&:= 
\Hom_{A'_{n \alpha_i}}(P(i^{(n)}), -): \mod{A_{\alpha+n\alpha_i}} \rightarrow \mod{A_{\alpha}},
%\label{div2}
\end{align*}
where $A'_{n\alpha_i} := 1 \otimes A_{n\alpha_i} \subseteq A_{\alpha,n\alpha_i}$. Define 
\begin{align*}
e_i^{(n)} &:= \pr \circ (\theta_i^*)^{(n)} \circ \infl:
\mod{H_{\alpha+n\alpha_i}} \rightarrow \mod{H_\alpha},\\
f_i^{(n)} &:= q^{n^2 - n(\La_0-\alpha,\alpha_i)}\pr \circ \theta_i^{(n)} \circ \infl 
:\mod{H_\alpha} \rightarrow \mod{H_{\alpha+n\alpha_i}}. 
\end{align*}
By \cite[Lemma 4.8]{BKllt}, $e_i^n\cong [n]_q^!e_i^{(n)}$, $f_i^n\cong [n]_q^!f_i^{(n)}$, where $[n]_q^!:=[1]_q\dots[n]_q$,  the functors $e_i^{(n)}, f_i^{(n)}$ are exact, and send finite dimensional (resp. projective) modules to finite dimensional (resp. projective)  modules. Finally, we need the degree shift functors 
\begin{equation*}
k_i^{n}:\mod{H_\alpha} \rightarrow \mod{H_\alpha},\ M\mapsto q^{n(\La-\alpha,\alpha_i)}M \qquad(n\in\Z).
\end{equation*}

Consider the (locally unital) algebra $H:=\oplus_{\al\in Q_+} H_\al$, the categories 
$$\mod{H} := \bigoplus_{\alpha \in Q_+} \mod{H_\alpha}\quad \text{and}\quad \proj{H} := \bigoplus_{\alpha \in Q_+} \proj{H_\alpha},
$$ 
and the Grothendieck groups 
$[\mod{H}] = \bigoplus_{\alpha \in Q_+} [\mod{H_\alpha}]$ with $\Laurent$-basis $B$, and $[\proj{H}] = \bigoplus_{\alpha \in Q_+} [\proj{H_\alpha}]$. 
Let 
\begin{equation}\label{cartanpairing}
\langle.,.\rangle:\proj{H}\times\mod{H} \rightarrow \Laurent,
\quad 
\langle[P],[M]\rangle := \DIM\ \Hom_H(P,M),
\end{equation}
be the {\em Cartan pairing}. 
The pairing is {\em sesquilinear}, i.e. anti-linear in the first argument and linear in the second. We have a similar  form $\langle.,.\rangle:\proj{H}\times\proj{H} \rightarrow \Laurent$. 

The exact functors $e_i^{(n)}, f_i^{(n)}$ and $k_i^{\pm 1}$ induce
$\Laurent$-linear endomorphisms $E_i^{(n)}, F_i^{(n)}$ and $K_i^{\pm 1}$, respectively of $[\mod{H}]$ and $[\proj{H}]$. 

\subsection{LLT categorification}
On the other hand, let $U_q(\g)$ be the quantized enveloping algebra of $\g$ over $\Q(q)$ with Chevalley generators $E_i^{(n)}, F_i^{(n)},K_i^{\pm1}$ for $i\in I$. 
Let $V(\La_0)$ be the irreducible $U_q(\g)$-module with highest weight $\La_0$ and a fixed highest weight vector $v_+$. The module $V(\La_0)$ has a unique
compatible bar-involution
$-:V(\La_0) \rightarrow V(\La_0)$
such that $\overline{v_+} = v_+$.

The {\em Shapovalov form}
$\langle.,.\rangle$ is the unique sesquilinear $\Q(q)$-valued form on $V(\La_0)$ such that $\langle v_\La, v_\La\rangle = 1$
and $\langle u v, w\rangle = \langle v, \tau(u) w\rangle$ for
all $u \in U_q(\g)$ and $v, w \in V(\La_0)$, where $\tau$ is 
anti-linear anti-automorphism defined by 
$\tau(K_i)= K_i^{-1}, \  \tau(E_i) = q F_i K_i^{-1},\ \tau(F_i) = q^{-1} K_i E_i$. 

Let $U_q(\g)_{\Laurent}$ be the Lusztig's $\Laurent$-form, i.e. 
the $\Laurent$-subalgebra of $U_q(\g)$ generated by the quantum divided
powers $E_i^{(n)} := E_i^n / [n]_q^!$ and $F_i^{(n)} := F_i^n / [n]_q^!$ for all $i \in I$ and $n \geq 1$. 
Let $V(\La_0)_{\Laurent}:=U_q(\g)_{\Laurent}\cdot v_+$ be the standard $\Laurent$-form of $V(\La_0)$, and $V(\La_0)_{\Laurent}^*= \{v \in V(\La_0)\:|\:\langle v,w \rangle
\in {\Laurent}\text{ for all }
w\in V(\La_0)_{\Laurent}\}$ be the costandard ${\Laurent}$-form. %} of $V(\La_0)$. 

The following is the graded version \cite[Theorem 4.18]{BKllt} of the categorification theorems proved by Lascoux-Leclerc-Thibon \cite{LLT}, Ariki \cite{Ariki}, and Grojnowsky \cite{Gr}:

\begin{Theorem} %\label{}%{\rm \cite{}}%{\bf ()}
The linear operators $E_i, F_i$ and $K_i$ on the Grothendieck group 
$$[\proj{H}]_{\Q(q)}:=[\proj{H}]\otimes_\Laurent\Q(q)$$ satisfy the defining relations of the Chevalley generators of $U_q(\g)$. So $[\proj{H}]_{\Q(q)}$ is a $U_q(\g)$-module. Moreover:
\begin{enumerate}
\item[{\rm (i)}] There is a unique isomorphism $\de:V(\La_0)\iso [\proj{H}]_{\Q(q)}$ of $U_q(\g)$-modules, such that $\de(v_+)=[\triv_{H_0}]$, where $\triv_{H_0}\in \proj{H_0}$ is the one-dimensional vector space $F$ considered as a module over $H_0\cong F$. 
\item[{\rm (ii)}] The restriction of $\de$ to $V(\La_0)_\Laurent$ is an  isomorphism $\de:V(\La_0)_\Laurent\iso[\proj{H}]$ of $U_q(\g)_\Laurent$-modules, which intertwines $\circledast$ with
the bar-involution
on $V(\La_0)_\Laurent$ and induces the isomorphisms on weight spaces $V(\La_0)_{\La_0-\al,\Laurent}\iso[\proj{H_\al}]$ for all $\al\in Q_+$. 
\item[{\rm (iii)}] The isomorphism $\de$ identifies 
the Shapovalov form on $V(\La_0)_\Laurent$
with the Cartan pairing on $[\proj{H}]$. 
\item[{\rm (iv)}] Let $\de^*:[\mod{H}]\to V(\La_0)_\Laurent^*$ be the dual map:
$$
\de^*([M])(v):=\langle \de(v),[M]\rangle\qquad(v\in V(\La_0)_\Laurent).
$$
Then $\de^*$ is an  isomorphism of $U_q(\g)_\Laurent$-modules, which intertwines $\circledast$ with
the bar-involution
on $V(\La_0)^*_\Laurent$, and induces the isomorphisms  $[\mod{H_\al}]\iso V(\La_0)_{\La_0-\al,\Laurent}^*$ for all $\al\in Q_+$. 

\item[{\rm (v)}] The following diagram is commutative:
$$
\begin{CD}
V(\La_0)_\Laurent&@>\sim > \de>&[\proj{H}]\\
@Va VV&&@VVbV\\
V(\La_0)_\Laurent^*&@<\sim <\de^*<&[\mod{H}],
\end{CD}
$$
where $a:V(\La)_\Laurent \hookrightarrow V(\La)_\Laurent^*$
is the canonical inclusion,
and
$b:[\proj{H}]  \rightarrow [\mod{H}]$
is the $\Laurent$-linear map induced by the natural inclusion
of $\proj{H}$ into $\mod{H}$. 
In particular, $b$ is injective and becomes an isomorphism
over $\Q(q)$.
\end{enumerate}
\end{Theorem}

We complete this section with a special case of the Chuang-Rouquier result \cite{CR} on derived equivalence of the algebras $H_\al$. Recall from \cite{Kac} that the (affine) Weyl group $W$ of $\g$ acts on the weights of $V(\La_0)$.

\begin{Theorem} %\label{}%{\rm \cite{}}%{\bf ()}
Let $\al,\be\in Q_+$. Then the derived categories $D^b(\mod{H_\al})$ and $D^b(\mod{H_\be})$ are equivalent if and only if the weights $\La_0-\al$ and $\La_0-\be$ belong to the same $W$-orbit. %, where $W$ is the (affine)  Weyl group of $\g$. 
\end{Theorem}

The equivalence in the theorem is induced by a complex of functors, which is built out of the functors $e_i$ and $f_i$ using adjunctions, see \cite[\S6]{CR}.

\section{Combinatorics of partitions and homogeneous representations}\label{SNot}

\subsection{Partitions and nodes}\label{SSPartNode}
Let $\Par_n$ be the set of all partitions of $n$  and put $\Par:=\bigsqcup_{n\geq 0}\Par_n$. If $\mu\in\Par_n$, we write $n=|\mu|$. A partition $\mu=(\mu_1,\mu_2,\dots)$ is called {\em $p$-restricted} if $\mu_k-\mu_{k+1}<p$ for all $k=1,2,\dots$. Let $\RPar_n$ be the set of all $p$-restricted partitions of $n$,  and put $\RPar:=\bigsqcup_{n\geq 0}\RPar_n$. 
%Let $\mu,\nu\in \Par_n$. Then $\mu$ {\em dominates} $\nu$, written $\mu\unrhd\nu$, if $\sum_{b=1}^c\mu_b\geq \sum_{b=1}^c\nu_b$ for all  $c\geq 1$.
The {\em Young diagram} of a partition $\mu=(\mu_1,\mu_2,\dots)$ is 
$
\{(a,b)\in\Z_{>0}\times\Z_{>0}\mid 1\leq b\leq \mu_a\}.
$ The elements of this set
are the {\em nodes of~$\mu$}. More generally, a {\em node} is any element of $\Z_{>0}\times\Z_{>0}$. 
We identify partitions with their 
Young diagrams, so that a node $(a,b)=$ box in row $a$ and column $b$. For example, %$(3,2,2,1)$ is 
$$
\begin{picture}(0,62)
\put(-55,25){$(3,2,2,1)=$}
\put(0,0){\YoungDiagram{3,2,2,1}}
\put(55,25){.}
\end{picture}
$$
To each node $A=(a,b)$ %of the Young diagram $\mu$ 
we associate its {\em residue}: 
$$
\res A:=(b-a)\pmod{p} \in I.
$$ 
An {\em $i$-node} is a node of residue $i$. Let $c_i(\mu)$ be the number of $i$-nodes of $\mu$, and define the {\em content} of $\mu$ to be 
$
\cont(\mu)=\sum_{i\in I} c_i(\mu)\al_i\in Q_+. 
$
Denote
$$\Par_\al:=\{\mu\in\Par\mid \cont(\mu)=\al\}, \quad \RPar_\alpha := \RPar \cap \Par_\alpha 
\qquad(\al\in Q_+).$$

%In other words, $\mu$ is obtained from $\nu$ by moving nodes up in the diagram.

A node $A\in\mu$ is {\em removable (for~$\mu$)}\, if $\mu\setminus \{A\}$ is a partition. A node $B\not\in\mu$ is an {\em addable node (for~$\mu$)}\, if $\mu\cup \{B\}$ is a partition. We denote 
$
\mu_A:=\mu\setminus \{A\},\ \mu^B:=\mu\cup\{B\}.
$

Let $i \in I$, and $A_1,\dots,A_n$ be the addable and removable $i$-nodes
of $\mu$ ordered so that $A_m$ is to the left of $A_{m+1}$ for each
$m=1,\dots,n-1$. 
Consider the sequence $(\tau_1,\dots,\tau_n)$
where $\tau_r = +$ if $A_r$ is addable and $-$ if $A_r$ is removable. 
If there are $1 \leq r < s \leq n$ with $\tau_r=  +$,
$\tau_s = -$ and $\tau_{r+1}=\cdots=\tau_{s-1} = 0$
then replace $\tau_r$ and $\tau_s$ by $0$.
Keep doing this until %we are 
left with a sequence
$(\sigma_1,\dots,\sigma_n)$ in which no $+$ appears to the left of
a $-$. This is the {\em reduced $i$-signature} of $\mu$ (it is well-defined). 

If $(\sigma_1,\dots,\sigma_n)$ is the reduced $i$-signature of $\mu$,
we set 
$$\eps_i(\mu) := \#\{r=1,\dots,n\:|\:\sigma_r = -\},\quad \phi_i(\mu) := \#\{r=1,\dots,n\:|\:\sigma_r = +\}.
$$ 
Let $\{r_1>\dots>r_{\eps_i(\mu)}\}=\{r\mid \si_r=-\}$, and $\{a_1<\dots<a_{\phi_i(\mu)}\}=\{a\mid \si_a=+\}$. 
If $\eps_i(\mu) > 0$, set 
$\tilde e_i \mu := \mu_{A_{r_1}}$; otherwise set $\tilde e_i \mu:=0$. If $\phi_i(\mu) >0 $, set 
$\tilde f_i \mu := \mu^{A_{a_1}}$; otherwise set $\tilde f_i \mu:=0$.  The removable nodes $A_{r_1},\dots,A_{r_{\eps_i(\mu)}}$ of $\mu$ are called {\em $i$-normal}, and the addable nodes $A_{a_1},\dots,A_{a_{\phi_i(\mu)}}$ of $\mu$ are called {\em $i$-conormal}. 

\subsection{Tableaux}
Let $\mu\in\Par_n$. 
A {\em $\mu$-tableau} 
$\T$ is %obtained from the diagram of $\mu$ by 
an insertion of the integers $1,\dots,n$ into the nodes of $\mu$, allowing no repeats. 
%If the node $A=(a,b)\in\mu$ is occupied by the integer $r$ in $\T$ then we write $r=\T(a,b)$ and set $\res_\T(r)=\res A$. 
The {\em residue sequence} of~$\T$ is
\begin{equation*}\label{EResSeq}
\bi^\T=(i_1,\dots,i_n)\in I^n,
\end{equation*}
where $i_r$ is the residue of the node occupied by 
$r$ in $\T$ ($1\leq r\leq n$). 
%\Comment[Andrew]{I added $\bi^\kappa(\T)$ here and similar changes to Lemma~8.4 below.}
A $\mu$-tableau $\T$ is {\em row-strict} (resp. {\em column-strict}) if its entries increase from left to right (resp. from top to bottom) along the rows (resp. columns) of $\mu$. 
A $\mu$-tableau $\T$ is {\em standard} if it is row- and column-strict. 
 Let $\St(\mu)$ be the set of all standard $\mu$-tableaux.

Let $\T^\mu$ be the {\em leading}\, $\mu$-tableau, i.e. the tableau in which the numbers $1,2,\dots,n$ appear in order from left to right along the successive rows,
working from top row to bottom row. For example, if $\mu=(3,2,2,1)$ then $\T^\mu$ is
$$
\Tableau{{1,2,3},{4,5},{6,7},{8}}.
$$
Set 
$%\begin{equation}\label{EBIMu}
\bi^\mu:=\bi^{\T^\mu}.
% \quad\text{and}\quad\bi_\mu:=\bi(\T_\mu).
$ %\end{equation} 
The group $\Si_n$ acts on the set of $\mu$-tableaux %on the left 
by acting on the entries of the tableaux.
For each $\mu$-tableau $\T$, define $w^\T\in\Si_n$ from  
$
w^\T  \T^\mu=\T
$. %\end{equation}

%Let $\ell$ be the length function on $\Si_d$ with respect to the Coxeter generators $s_1,s_2,\dots,s_{d-1}$. 
Let $\le$ be the {\em Bruhat order} on $\Si_d$. Define the {\em Bruhat order} on the set of all $\mu$-tableaux as follows: %if $\Stab,\T\in\St(\mu)$ then
$ %\begin{equation}\label{EBruhat}
\Stab\ledom\T\quad\text{if and only if}\quad w^\Stab\le w^\T.
$ %\end{equation}
%If $\Stab\ledom\T$ then we also write $\T\gedom\Stab$. If $\Stab\ledom\T$ and $\Stab\ne\T$ we write $\Stab\ldom\T$ and $\T\gdom\Stab$.
Then the leading $\mu$-tableau $\T^\mu$ is the unique minimal element of $\St(\mu)$. 

Let $\mu\in\Par$, $i\in I$, and $A$ be a removable $i$-node %and $B$ be an addable $i$-node 
of $\mu$. We set
\begin{align*}\label{EDMUA}
d_A(\mu)&=\#\SetBox{addable $i$-nodes of $\mu$\\[-4pt] strictly to the left of $A$}
                -\#\SetBox[38mm]{removable $i$-nodes of $\mu$\\[-4pt] strictly to the left of $A$}.
\end{align*} 
Given $\T \in \St(\mu)$, the {\em degree} of $\T$ is defined in
\cite[\S3.5]{BKW} inductively as follows. If $n=0$, %then $\T$ is the empty tableau $\emptyset$, and
we set $\deg(\T):=0$.  Otherwise, let $A$ be the node occupied by $n$ in~$\T$. Let $\T_{<n}\in\St(\mu_A)$
be the tableau obtained by removing $A$ and set
$$%\begin{equation}\label{EDegTab}
\deg(\T):=d_A(\mu)+\deg(\T_{<n}).
$$ %\end{equation}

\section{Branching and graded cellular structure}

\subsection{Crystal combinatorics and irreducible modules}
Using the terminology introduced in Section~\ref{SNot}, we can now state the following theorem of Misra and Miwa \cite{MisMiw}: 

\begin{Theorem} \label{TCrystComb}%{\rm \cite{}}%{\bf ()}
For any partition $\mu$, define $\wt(\mu):=\La_0-\cont(\mu)$. Then the tuple $(\RPar,\eps_i, \phi_i,\tilde e_i, \tilde f_i,  \wt)
$ is the Kashiwara's crystal associated to $V(\La_0)$.
\end{Theorem}

Comparing this with Theorem~\ref{TCryst}, we deduce that there is a unique isomorphism of crystals 
\begin{equation}\label{CrystIso}
(\RPar,\eps_i, \phi_i,\tilde e_i, \tilde f_i,  \wt)
\iso (B,\eps_i,\phi_i,\tilde e_i, \tilde f_i,\wt). 
\end{equation}
Under this isomorphism, to every $\mu\in\RPar_\al$, we associate the irreducible $H_\al$-module $D^\mu$, and 
\begin{equation}\label{EDMu}
\{D^\mu\mid \mu\in\RPar_\al\}
\end{equation} 
is a complete and irredundant set of irreducible $H_\al$-modules up to isomorphism and degree shift. %Moreover, we can now reformulate 

On the other hand, there is another approach to the classification of irreducible $H_\al$-modules, based on the theory of {\em Specht modules}, and which goes back to James \cite{JamesB}. In modern terms, this is the approach through {\em cell modules}\, \cite{GrLeh}. The {\em graded cellular structure}\, of $H_\al$, which we present here, has been worked out by Hu and Mathas \cite{HM}. 

\subsection{Graded cellular structure}
Fix $\al\in Q_+$ and $\mu\in\Par_\al$. Recall the leading standard tableau $\T^\mu\in\St(\mu)$ and the corresponding residue sequence $\bi^\mu$. For $k=1,\dots,n$, let $A_k$ be the box occupied with $k$ in $\T^\mu$. Observe that $A_k$ is a removable node for the partition $\mu^k$, obtained from $\mu$ by removing $A_{k+1},\dots, A_n$. Set $d_k(\mu):=d_{A_k}(\mu^k)$. Note that $\deg(T^\mu)=d_1(\mu)+\dots+d_n(\mu)$. Define
$$
y^\mu:=y_1^{d_1(\mu)}\dots y_n^{d_n(\mu)}.
$$

Given $\T\in\St(\mu)$, recall the element $w^\T\in\Si_n$. Pick any reduced decomposition $w^\T=s_{m_1}\dots s_{m_l}$, and set 
\begin{equation}\label{EPsiT}
\psi^{\T}:=\psi_{m_1}\dots\psi_{m_l}\in H_\al. 
\end{equation}
This element in general depends on the choice of the reduced decomposition. Finally, recalling (\ref{EPrime}), for any standard tableaux $\Stab,\T\in\St(\mu)$, we define
$$
\psi^{\Stab,\T}:=\psi^{\Stab}y^\mu 1_{\bi^\mu} (\psi^{\T})'.
$$
It is easy to see that $\deg(\psi^{\Stab,\T})=\deg(\Stab)+\deg(\T)$. 
The following theorem was proved by Hu and Mathas \cite{HM}:

\begin{Theorem} %\label{}%{\rm \cite{}}%{\bf ()}
Let $\al\in Q_+$. Then $\{\psi^{\Stab,\T}\mid \mu\in\Par_\al,\ \Stab,\T\in\St(\mu)\}$ is a graded cellular basis of $H_\al$. 
\end{Theorem}

The following immediate corollary was originally proved in \cite[Theorem 4.20]{BKllt} by a different method:

\begin{Corollary} %\label{}%{\rm \cite{}}%{\bf ()}
Let $\al\in Q_+$ and $\bi,\bj\in\al$. Then 
$$\displaystyle \DIM 1_\bi H_\al1_\bj=\sum_{\mu\in\Par_\al,\ \Stab,\T\in\St(\mu),\ \bi^{\Stab}=\bi,\ \bj^\T=\bj}q^{\deg(\Stab)+\deg(\T)}.$$ 
In particular, 
$\displaystyle \DIM H_\al=\sum_{\mu\in\Par_\al,\ \Stab,\T\in\St(\mu)}q^{\deg(\Stab)+\deg(\T)}$. 
\end{Corollary}

The graded version of the Graham-Lehrer theory \cite{GrLeh} can be found in \cite{HM}. In particular, from a graded cellular basis we get graded cell modules $\{S^\mu\mid \mu\in\Par_\al\}$. It is shown in \cite{HM} that these are just the {\em graded Specht modules} as constructed originally in \cite{BKW}. On the other hand, it is noted in \cite{BKW} that, if we forget the grading, then $S^\mu$ is the usual dual Specht module of \cite{JamesB}. 

Recall from (\ref{EDMu}) the simple $H_\al$-modules $D^\mu$  defined using the crystal isomorphism~(\ref{CrystIso}).

\begin{Theorem} %\label{}%{\rm \cite{}}%{\bf ()}
If $\mu$ is $p$-restricted, then $S^\mu$ has a simple head isomorphic to $D^\mu$.  
\end{Theorem}

There might be more to this theorem than meets the eye. Firstly, it is the graded aspect: recall that $(D^\mu)^\circledast\cong D^\mu$, and the theorem claims that the natural map $S^\mu\surj \head S^\mu\cong D^\mu$ is a homogeneous degree zero map. Secondly, it is known from James \cite{JamesB} that the head of $S^\mu$ is simple when $\mu$ is $p$-restricted, and this is more or less how James classifies the irreducible modules (actually our $D^\la$ is isomorphic to James' $D^{\la'}\otimes \sgn$). But it is not at all clear that the James classification agrees with the classification (\ref{EDMu}) coming from the isomorphism of crystals (\ref{CrystIso}). There are several ways of seeing this, none being trivial. One comes from the fact that the original branching rules are proved in \cite{KBrII} for the modules $D^\la$ in James' classification, which allows us to identify them with the classification (\ref{EDMu}). 

\subsection{Mullineux map}
Branching rules yield a simple solution to the Mullineux problem \cite{Mull}. Tensoring with the sign representation yields a bijection $$\RPar(\al)\iso \RPar(-\al),\ \mu\mapsto\Mull(\mu),$$ 
where $\Mull(\mu)$ is defined from $$D^{\Mull(\mu)}\cong D^\mu\otimes\sgn\qquad(\mu\in\RPar_\al).$$ The problem is to describe $\Mull(\mu)$ explicitly in combinatorial terms. The following theorem, proved in \cite{KBrIII}, gives an answer to this question:

\begin{Theorem} %\label{}%{\rm \cite{}}%{\bf ()}
Let $\mu\in\RPar_n$, and $\emptyset$ be the empty partition. Pick a sequence $i_1,\dots,i_n\in I$ such that $\emptyset=\tilde e_{i_1}\dots \tilde e_{i_n}\mu$. Then $\Mull(\mu)=\tilde f_{-i_n}\dots\tilde f_{-i_1}\emptyset$.
\end{Theorem}

Thus a computation of $\Mull(\mu)$ is reduced to combinatorics of the crystal graph of Theorem~\ref{TCrystComb} described in \S\ref{SSPartNode}. However, there a faster algorithm originally conjectured by Mullineux \cite{Mull}. The Mullineux Conjecture has been first proved in \cite{FK} and simpler proofs were later found in \cite{BO} and \cite{BrKu}. We now describe this algorithm or rather its more elegant  version suggested by Xu \cite{Xu}. 

The {\em rim} of $\mu$ is the set of all nodes $(i, j) \in\mu$ such that that $(i + 1, j + 1) \not\in\mu$. The {\em $p$-rim} of $\mu$ is the  union of the {\em $p$-segments}, which are defined as follows. The first $p$-segment of $\mu$ consists of the first $p$ nodes of the rim, reading along the rim from bottom-left to top-right. The next $p$-segment is obtained by similar reading off the nearest $p$ nodes of the rim, but {\em starting from the column immediately to the right of the rightmost node of the first $p$-segment}. And so on. All but the last $p$-segment contain exactly $p$ nodes, while the last may contain less.
In the following example $p=3$, there are two $p$-segments, and the nodes of the $p$-rim are marked with $*$'s. 
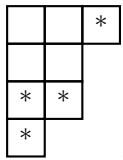
$$
\Tableau{{,,*},{,},{*,*},{*}}.
$$

Let let $J(\mu)$ be the partition obtained from $\mu$ by deleting every node in the $p$-rim that is at the rightmost end of a row of $\mu$ but that is not the $p$th node of a $p$-segment. Let $j(\mu) = |\mu| - |J(\mu)|$ be the total number of nodes deleted. Now, in Xu's reformulation \cite{Xu}, the result is as follows:

\begin{Theorem} %\label{}%{\rm \cite{}}%{\bf ()}
$\Mull(\mu)$ is the partition $\la=(\la_1,\la_2,\dots)$ with $\la_r = j(J^{r-1}(\mu))$. 
\end{Theorem}

In the example above, we get $\Mull((3,2^2,1))=(2,1^6)$.

\subsection{More branching rules}\label{SSMoreBR}
We complete this section with more result on branching. 

\begin{Theorem} \label{TNorm}%{\rm \cite{}}%{\bf ()}
Let $\al\in Q_+$, $i\in I$, $\mu\in\RPar_\al$, $A$ be a removable node of $\mu$ such that $\mu_A$ is $p$-restricted, and $B$ be an addable node for $\mu$ such that $\mu^B$ is $p$-restricted. Moreover, let $A_1,\dots,A_{\eps_i(\mu)}$ be the $i$-normal nodes of $\mu$ counted from left to right, and $B_1,\dots,B_{\phi_i(\mu)}$ be the $i$-conormal nodes for $\mu$ counted from right to left. 
\begin{enumerate}
\item[{\rm (i)}] $\Hom_{H_{\al-\al_i}}(S^{\mu_A},e_i D^\mu)\neq 0$ if and only if $A$ is $i$-normal for $\mu$, in which case we have $\DIM \Hom _{H_{\al-\al_i}}(S^{\mu_{A_m}},e_i D^\mu)=q^{m-1}$ for all $m=1,\dots,\eps_i(\mu)$. % such that $\mu_{A_m}$ is $p$-restricted. 

\item[{\rm (ii)}] $\Hom_{H_{\al+\al_i}}(S^{\mu^B},f_i D^\mu)\neq 0$ if and only if $B$ is $i$-conormal for $\mu$, in which case we have $\DIM \Hom _{H_{\al+\al_i}}(S^{\mu^{B_m}},f_i D^\mu)=q^{m-1}$ for all $m=1,\dots,\phi_i(\mu)$. % such that $\mu^{B_m}$ is $p$-restricted. 

\item[{\rm (iii)}] $D^{\mu_A}$ appears as a composition factor of $e_i D^\mu$ if and only if $A$ is $i$-normal for $\mu$, in which case we have $[e_i D^\mu:D^{\mu_{A_m}}]_q=[m]_q$  for all $m=1,\dots,\eps_i(\mu)$. % such that $\mu_{A_m}$ is $p$-restricted. 

\item[{\rm (iv)}] $D^{\mu^B}$ appears as a composition factor of $f_i D^\mu$ if and only if $B$ is $i$-conormal for $\mu$, in which case we have $[f_i D^\mu:D^{\mu^{B_m}}]_q=[m]_q$  for all $m=1,\dots,\phi_i(\mu)$. % such that $\mu^{B_m}$ is $p$-restricted. 
\end{enumerate}
\end{Theorem}

This is a graded version of the results \cite{KBrII}, \cite{KDec},\cite{KBrIV},  \cite{BKTr}. The graded version is deduced using the graded endomorphism algebra description of Theorem~\ref{TBr}(v) and the related filtrations of $e_i D^\mu$ and $f_i D^\mu$ obtained in \cite{KBrIV} and \cite{BKTr}. The following corollary follows immediately from Theorem~\ref{TNorm} on tensoring with sign, cf. \cite{KDec}, and often provides us with some new non-trivial branching multiplicities:

\begin{Corollary} \label{CNorm}%{\rm \cite{}}%{\bf ()}
Let $\al\in Q_+$, $i\in I$, $\mu\in\RPar_\al$, and let $A_1,\dots,A_{\eps_i(\mu)}$ be the $(-i)$-normal nodes of $\Mull(\mu)$ labeled from left to right, and $B_1,\dots,B_{\phi_i(\mu)}$ be the $(-i)$-conormal nodes for $\Mull(\mu)$ labeled from right to left. Then: 
\begin{enumerate}
\item[{\rm (i)}] $[e_i D^\mu:D^{\Mull(\Mull(\mu)_{A_m})}]_q=[m]_q$  for all $m=1,\dots,\eps_i(\mu)$ such that $\Mull(\mu)_{A_m}$ is $p$-restricted. 

\item[{\rm (ii)}] $[f_i D^\mu:D^{\Mull(\Mull(\mu)^{B_m})}]_q=[m]_q$  for all $m=1,\dots,\phi_i(\mu)$ such that $\Mull(\mu)^{B_m}$ is $p$-restricted. 
\end{enumerate}
\end{Corollary}

Let us consider the example where $p=3$ and $\mu=(3^2,2,1^2)$. We draw the corresponding Young diagram with the residues of the boxes written in them.   
$$
\Tableau{{0,1,2},{2,0,1},{1,2},{0},{2}}.
$$
Note that both $2$-removable boxes are normal, the removable $1$-box is not, and there are no removable $0$-boxes. So $e_0 D^\mu=e_1 D^\mu=0$. As for the composition factors of $e_2 D^\mu$, Theorem~\ref{TNorm} shows that $[e_2 D^\mu:D^{(3^2,2,1)}]_q=1$ and $[e_2 D^\mu:D^{(3^2,1^3)}]_q=[2]_q=q+q^{-1}$. Moreover, since $\Mull(\mu)=(3^2,1^4)$ has the leftmost normal $1$-node $(6,1)$ and $\Mull(3^2,1^3)=(3,2,1^4)$, Corollary~\ref{CNorm} yields another composition factor $D^{(3,2,1^4)}$ of multiplicity $1$. It is easy to verify using decomposition matrices in \cite{JamesB} that in this example %, using only Theorem~\ref{TNorm} and Corollary~\ref{CNorm}, 
we have discovered all composition factors, i.e. 
$[\res^{\Si_{10}}_{\Si_{9}}]= [D^{(3^2,2,1)}]+(q+q^{-1})[D^{(3^2,1^3)}]+[D^{(3,2,1^4)}].$

 Unfortunately, this technique is not powerful enough to always yield all composition factors, see \cite[Section 1]{KDec}. So it leads only to a lower bound on the dimensions of  irreducible $H_n$-modules.  A  family of irreducible modules for which this lower bound is equal to the actual dimension is described in the next section. 

\subsection{Homogeneous representations}\label{SSHomog}
Let $\al\in Q_+$. An irreducible $H_\al$-module is called {\em homogeneous} if it is concentrated in degree zero. To describe the homogeneous representations, for a partition $\mu=(\mu_1\geq.\dots\geq \mu_u>0)\in\RPar_n$ consider the hook length 
$\chi(\mu):=\la_1+u-\max\{t\mid \la_t=\la_1\}$. 
Then $\mu$ is called {\em homogeneous} if $\chi(\mu)\leq p$, cf. \cite{KCS}, where we worked with transposed partitions. 

Let $\mu\in\RPar_n$ be a homogeneous partition. A tableaux $\T\in\St(\mu)$ is called {\em $p$-standard} if 
%$\T_{<m}$ is homogeneous for all $m=1,\dots,n$ (as usual 
$a<b$ whenever $a$ occupies a box $(r,s)$ in $\T$ and $b$ occupies a box $(r',s')$ with $r>r'$, $s<s'$ and $r-r'+s'-s+1=p$. Let $\St^p(\mu)$ be the set of all $p$-standard $\mu$-tableaux. The results of \cite{KCS,Mat} and \cite{KRhomog} can be restated as follows:

\begin{Theorem} \label{Thomog}%{\rm \cite{}}%{\bf ()}
Let $\mu\in \RPar_\al$. 
The irreducible $H_\al$-module $D^\mu$ is homogeneous if and only if $\mu$ is a homogeneous partition. In this case, $D^\mu$ has a basis $\{v_\T\mid \T\in\St^p(\mu)\}$ with the action of the homogeneous  generators of $H_\al$ given as follows:
$$
1_\bi v_\T=\de_{\bi,\bi^\T}v_\T,\ y_t v_\T=0, \ \psi_rv_\T=
\left\{
\begin{array}{ll}
v_{s_r \T} &\hbox{if $s_r \T\in \St^p(\mu)$;}\\
0 &\hbox{otherwise.}
\end{array}
\right.
$$
\end{Theorem}

\section{Presentations and bases of cell modules}

\subsection{Garnir tableaux}
Let  $\mu\in\Par_n$. We now explain an explicit presentation of the graded Specht module $S^\mu$ obtained in \cite{KMR}. First, we need more notation. Let $A=(r,s)$ be a node of $\mu$. It is called  a {\em  Garnir node} if $(r+1,s)\in\mu$, i.e. $A$ is not at the bottom of its column. Then the {\em $A$-Garnir belt}\, $\Belt^A$ is $$\Belt^A:=\set{(r,t)\in\mu| s\leq t\leq \mu_{r}}\cup
          \set{(r+1,u)\in\mu|1\leq u\leq s}. $$
For example,
if $A=(2,3)$ then $\Belt^A$ for $\mu=(7,7,4,1)$ is highlighted below:  
$$\begin{array}{l}
%\YoungDiagram{1}\\[10pt]
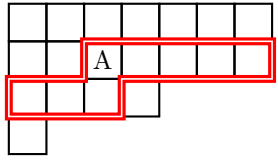
\begin{tikzpicture}[scale=0.5,draw/.append style={thick,black}]
  \newcount\col
  \foreach\Row/\row in {{,,,,,,}/0,{,,A,,,,}/-1,{,,,}/-2,{{}}/-3} {
     \col=1
     \foreach\k in \Row {
        \draw(\the\col,\row)+(-.5,-.5)rectangle++(.5,.5);
        \draw(\the\col,\row)node{\k};
        \global\advance\col by 1
      }
   }
   \draw[red,double,very thick]
     (2.5,-0.5)--++(5,0)--++(0,-1)--++(-4,0)--++(0,-1)--++(-3,0)--++(0,1)--++(2,0)--cycle;
\end{tikzpicture}
\end{array}$$
The {\em $A$-Garnir tableau}  is the $\mu$-tableaux $\G^A$ defined as follows. Let $u=\T^\mu(r,s)$ be the entry of the leading $\mu$-tableau $\T^\mu$ which occupies the node $A=(r,s)$, and $v=\T^\mu(r+1,s)$. To get $\G^A$, insert  the numbers $u,u+1,\dots,v$ into the nodes of the Garnir belt going from left bottom to top right, and the other numbers into the same positions as in $\T^\mu$. 
Continuing the previous example, %For example, if $\mu=((1),(7,7,4,1))$ then 
$u=10, v=17$, and: % $T^\mu$ and $\G^A$ are: 
$$
\T^\mu=\begin{array}{l}
%\Tableau{{1}}\\[10pt]
\begin{tikzpicture}[scale=0.5,draw/.append style={thick,black}]
  \newcount\col
  \foreach\Row/\row in {{1,2,3,4,5,6,7}/0,{8,9,10,11,12,13,14}/-1,{15,16,17,18}/-2,{19}/-3} {
     \col=1
     \foreach\k in \Row {
        \draw(\the\col,\row)+(-.5,-.5)rectangle++(.5,.5);
        \draw(\the\col,\row)node{\k};
        \global\advance\col by 1
      }
   }
   \draw[red,double,very thick]
     (2.5,-0.5)--++(5,0)--++(0,-1)--++(-4,0)--++(0,-1)--++(-3,0)--++(0,1)--++(2,0)--cycle;
\end{tikzpicture}
\end{array},\quad
\G^A=\begin{array}{l}
%\Tableau{{1}}\\[10pt]
\begin{tikzpicture}[scale=0.5,draw/.append style={thick,black}]
  \newcount\col
  \foreach\Row/\row in {{1,2,3,4,5,6,7}/0,{8,9,13,14,15,16,17}/-1,{10,11,12,18}/-2,{19}/-3} {
     \col=1
     \foreach\k in \Row {
        \draw(\the\col,\row)+(-.5,-.5)rectangle++(.5,.5);
        \draw(\the\col,\row)node{\k};
        \global\advance\col by 1
      }
   }
   \draw[red,double,very thick]
     (2.5,-0.5)--++(5,0)--++(0,-1)--++(-4,0)--++(0,-1)--++(-3,0)--++(0,1)--++(2,0)--cycle;
\end{tikzpicture}
\end{array}$$
%We let $\bi^A:=\bi^{\G^A}$. 

Fix a Garnir node $A=(r,s)$ of $\mu$. 
A {\em brick} is a set of $p$ successive nodes in the same row 
$\{(t,u),(t,u+1),\dots,(t,u+p-1)\}\subseteq \Belt^A$ such that
$\res(t,u)=\res A$. Note that $\Belt^A$ is a disjoint union of the bricks that it contains, together
with less than $p$ nodes at the end of row $r$ which are not contained in a brick and less than $p$ nodes
at the beginning of row $r+1$ which are not contained in a brick. 

Let $k%=k^A
$ be the number of bricks in $\Belt^A$ (possibly zero). We label the bricks 
$B_1^A,B_2^A,\dots, B_k^A$ 
%be the bricks contained in the Garnir belt $\Belt$, ordered
going from left to right along row $r+1$ and then from left to right along row $r$ of $\G^A$. % as in the example above. %Of course, it might happen that $\Belt^A$ does not contain any bricks (this is always true if $e=0$), in which case $k=0$.
For example, the following picture shows the bricks in the $(2,3)$-Garnir belt of
$\mu=(7,7,4,1)$ when $p=2$: 
$$\begin{array}{l}
%\Tableau{{1}}\\[-20pt]
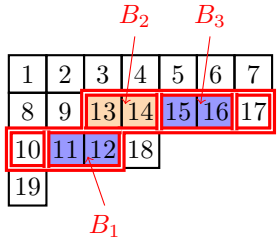
\begin{tikzpicture}[scale=0.5,draw/.append style={thick,black}]
  \fill[orange!30](2.5,-1.5)rectangle(4.5,-0.5);
  \fill[blue!40](4.5,-1.5)rectangle(6.5,-0.5);
  \fill[blue!40](1.5,-2.5)rectangle(3.5,-1.5);
  \newcount\col
  \foreach\Row/\row in {{1,2,3,4,5,6,7}/0,{8,9,13,14,15,16,17}/-1,{10,11,12,18}/-2,{19}/-3} {
     \col=1
     \foreach\k in \Row {
        \draw(\the\col,\row)+(-.5,-.5)rectangle++(.5,.5);
        \draw(\the\col,\row)node{\k};
        \global\advance\col by 1
      }
   }
   \draw[red,double,very thick](2.5,-1.5)--(2.5,-0.5)--(7.5,-0.5)--(7.5,-1.5)--(3.5,-1.5)
                     --(3.5,-2.5)--(0.5,-2.5)--(0.5,-1.5)--(2.5,-1.5);
   \draw[red,double,very thick](4.5,-0.5)--(4.5,-1.5);
   \draw[red,double,very thick](6.5,-0.5)--(6.5,-1.5);
   \draw[red,double,very thick](1.5,-1.5)--(1.5,-2.5);
   \draw[red,double,very thick](2.5,-1.5)--(3.5,-1.5);
   \draw[red,thin,->](3.85,1)node[red,above=0mm]{$B_2$}--(3.6,-0.8);
   \draw[red,thin,->](5.85,1)node[red,above=0mm]{$B_3$}--(5.6,-0.8);
   \draw[red,thin,->](3.05,-3.5)node[red,below=0mm]{$B_1$}--(2.6,-2.2);
\end{tikzpicture}
\end{array}$$
We have $k=3$, there are two bricks $B_2$, $B_3$ in row $2$ and one brick $B_1$ in row $3$.  
Finally, $(3,1)$ and $(2,7)$ are the nodes in the Garnir belt which are not contained in a brick.

Assume now that $k>0$ and let $d$ be the smallest entry in $\G^A$ which is contained in a brick in
$\Belt^A$. In the example above, $n=11$. Define
\begin{equation}%\label{EWRA}
w_t^A=\prod_{a=d+tp-p}^{d+tp-1}(a,a+p)\in\Si_n\qquad (1\le t<k).
\end{equation}
Informally, $w_t^A$ is a brick permutation swapping 
$B_t$ and $B_{t+1}$. The elements 
$w_1^A,w_2^A,\dots,w_{k-1}^A$ are the Coxeter generators of the
{\em brick permutation group}  
$$\Si^A:=\langle w_1^A,w_2^A,\dots,w_{k-1}^A\rangle\cong \Si_k.$$ 
%We call $\Si^A$ the {\em brick permutation group}.  
By convention, $\Si^A$ is the trivial group if $k=0$. 

Let $\Gar^A$ be the set of all row-strict $\mu$-tableaux which are obtained from the Garnir tableau $\G^A$ by acting with the brick permutation group $\Si^A$ on $\G^A$. 
Note that all of the tableaux in $\Gar^A$, except for $\G^A$, are standard. Moreover, $\G^A$ is the maximal element of $\Gar^A$, with respect to the Bruhat order, and there is a unique minimal tableaux $\T^A$ in $\Gar^A$. 
%If $\T\in\Gar^A$ then $\bi^\T=\bi^{\G^A}$. Consequently, we let $\bi^A=\bi^{\G^A}$. % be this common residue sequence. 

Let $f$ be the number of bricks in row $r$ of the Garnir belt $\Belt^A$, and let $\D^A$ be
the set of minimal length left coset representations of $\Si_f\times\Si_{k-f}$ in $\Si^A\cong\Si_k$. By definition $\Si^A$ is a subgroup of $\Si_n$, so $\D^A$ is a subset of $\Si_n$, and, in particular, its elements act on $\mu$-tableaux. We have 
\begin{equation}\label{EGarD}
\Gar^A=\{w\T^A\mid w\in\D^A\}. 
\end{equation}

Continuing the example above, $\T^A$ is the tableau
$$\T^A=\begin{array}{l}
%\Tableau{{1}}\\[-20pt]
\begin{tikzpicture}[scale=0.5,draw/.append style={thick,black},baseline=1mm]
  \fill[orange!50](2.5,-1.5)rectangle(4.5,-0.5);
  \fill[blue!30](4.5,-1.5)rectangle(6.5,-0.5);
  \fill[blue!30](1.5,-2.5)rectangle(3.5,-1.5);
  \newcount\col
  \foreach\Row/\row in {{1,2,3,4,5,6,7}/0,{8,9,11,12,13,14,17}/-1,{10,15,16,18}/-2,{19}/-3} {
     \col=1
     \foreach\k in \Row {
        \draw(\the\col,\row)+(-.5,-.5)rectangle++(.5,.5);
        \draw(\the\col,\row)node{\k};
        \global\advance\col by 1
      }
   }
   \draw[red,double,very thick](2.5,-1.5)--(2.5,-0.5)--(7.5,-0.5)--(7.5,-1.5)--(3.5,-1.5)
                     --(3.5,-2.5)--(0.5,-2.5)--(0.5,-1.5)--(2.5,-1.5);
   \draw[red,double,very thick](4.5,-0.5)--(4.5,-1.5);
   \draw[red,double,very thick](6.5,-0.5)--(6.5,-1.5);
   \draw[red,double,very thick](1.5,-1.5)--(1.5,-2.5);
   \draw[red,double,very thick](2.5,-1.5)--(3.5,-1.5);
   %\draw[red,thin,->](3.85,1)node[red,above=0mm]{$B^A_1$}--(3.6,-0.8);
   %\draw[red,thin,->](5.85,1)node[red,above=0mm]{$B^A_2$}--(5.6,-0.8);
   %\draw[red,thin,->](3.05,-3.5)node[red,below=0mm]{$B^A_3$}--(2.6,-2.2);
\end{tikzpicture},
\end{array}$$
and 
$\Gar^A=\{\T^A,w_2^A\T^A,\G^A=w_1^Aw_2^A\T^A\}.$

\subsection{Presenting Specht modules}
Set 
%\begin{equation}\label{ESiTauA}
%\si_r^A:=\psi_{w_r^A}e(\bi^A)\quad\text{and}\quad 
$\tau_r^A:=(\psi_{w_r^A}+1).$ %e(\bi^A),
%\end{equation}
Any element $u\in\Si^A$ can be written as a reduced product $u=w_{r_1}^A\dots w_{r_a}^A$ of simple generators $w_1^A,\dots,w_{k-1}^A$ of $\Si^A$. 
%In general, the elements $\tau_r^A$ do not have to satisfy Coxeter relations. However, if $u$ is fully commutative then the element
Then define 
$
\tau_u^A:=\tau_{r_1}^A\dots \tau_{r_a}^A. 
$
%is well-defined, since $\tau_r^A$ and $\tau_s^A$ commute for $|r-s|>1$. 
%In particular, we have well-defined elements 
%$$\{\tau_u^A\mid u\in \D^A\}.$$
%({\em As  operators on the brick permutation space $T^{\mu,A}$}, defined below, the elements $\tau_r^A$ do satisfy Coxeter relations, see Theorem~\ref{TTauA}(ii).)   

%Recall from (\ref{EGarD}) that $\Gar^A$ is the set of row-strict tableaux obtained from the tableau $\T^A$ by acting with the elements of $\D^A$. Note that for any $\Stab\in\Gar^A$, we can write $w^\Stab=u^\Stab w^{\T^A}$ so that $\ell(w^\Stab)=\ell(u^\Stab)+\ell(w^{\T^A})$ and $u^\Stab\in\D^A$. Moreover, in view of Lemma~\ref{LFullComm}, all elements $w^\Stab,u^\Stab,$ and $w^{\T^A}$ are fully commutative so the elements $\psi_{u^\Stab}$, $\psi^{\T^A}$ and $\psi^\Stab=\psi_{u^\Stab}\psi^{\T^A}$ are all independent of the choice of preferred reduced decomposition. Set 
%$$m^A:=m^{\T^A}=\psi^{\T^A}m^\mu\in M^\mu. $$

%\begin{Definition}%\label{}%{\rm \cite{}}%{\bf ()}
Suppose that $\mu\in\Par_\al$ and $A\in\mu$ is a Garnir node. The {\em Garnir element} is  
$$
g^A:=%\sum_{u\in\D^A}\tau_u^Am^A=
\sum_{u\in\D^A}\tau_u^A\psi^{\T^A}1_{\bi^\mu}\in H_\al.%m^\mu\in M^\mu
$$
%\end{Definition}
%In the module $M^\mu$ we have  $$g^Am^\mu=\sum_{u\in\D^A}\tau_u^Am^A.$$ By Lemma~\ref{L10910} all of the summands on the right hand side have the same degree.  Finally, if $\D^A=\{1\}$, we have $\G^A=\T^A$ and $g^A=\psi^{\G^A}$. 
%This includes the case where there are no bricks. In particular, this is the case for all Garnir nodes
%$A\in\mu$ when $e=0$. 
A special case of the main result of \cite{KMR} is:

\begin{Theorem} \label{TKMR}%{\rm \cite{}}%{\bf ()}
Let $\al\in Q_+$ be of height $n$ and $\mu\in\Par_\al$. The graded Specht module $S^\mu$ has a homogeneous vector $z^\mu$ of degree $\deg(\T^\mu)$ such that $S^\mu$ is generated as a graded $H_\al$-module by $z^\mu$ subject only to the following relations:
\begin{enumerate}
\item[{\rm (i)}] $1_\bj z^\mu=\de_{\bj,\bi^\mu} z^\mu$ for all $\bj\in \al$;
\item[{\rm (ii)}] $y_rz^\mu=0$ for all $r=1,\dots,n$;
\item[{\rm (iii)}] $\psi_rz^\mu=0$ for all $1\leq r<n$ such that  $r$ and $r+1$ are in the same row of $\T^\mu$;
\item[{\rm (iv)}] {\em (homogeneous Garnir relations)}\, $g^A z^\mu=0$ for all Garnir nodes\, $A$ in~$\mu$. 
\end{enumerate}
\end{Theorem}

In other words, the theorem says that $S^\mu= q^{\deg(\T^\mu)} H_\al/J_\al^\mu$, where $J_\al^\mu$ is the homogeneous left ideal of $H_\al$ generated by the elements 
{\rm (i)} $1_\bj-\de_{\bj,\bi^\mu}$ for all $\bj\in I^\al$;
{\rm (ii)} $y_r$ for all $r=1,\dots,n$;
{\rm (iii)} $\psi_r$ for all $1\leq r<n$ such that  $r$ and $r+1$ are in the same row of $\T^\mu$;
{\rm (iv)} $g^A$\, for all Garnir nodes\, $A\in \mu$. 

We refer the reader to \cite{Fayers} for further developments on this presentation. 
A homogeneous basis of $S^\mu$ can now be given as follows \cite{BKW}:

\begin{Theorem} %\label{}%{\rm \cite{}}%{\bf ()}
Let $\al\in Q_+$, $\mu\in\Par_\al$ and $z^\mu\in S^\mu$ be the element from Theorem~\ref{TKMR}. For $\T\in\St(\mu)$, define $
v^\T:=\psi^\T z^\mu$.  
Then $v^\T$ is a homogeneous vector of degree  $\deg(v^\T)=\deg(\T)$, and 
$
\{v^\T\mid \T\in\St(\mu)\}
$
is a basis of $S^\mu$. In particular, $$\DIM S^\mu=\sum_{\T\in\St(\mu)}q^{\deg(\T)}.$$ 
\end{Theorem}

\section{Representation theory of KLR algebras}\label{SKLR}

We now return to the KLR algebra $A_\al$, defined in \S\ref{SSKLRDef}. Graded $H_\al$-modules inflate to graded $A_\al$-modules via the natural surjection (\ref{ESurj}). In particular, the irreducible modules $D^\mu$ inflate to irreducible $A_\al$-modules. However, representation theory of $A_\al$ is more rich and perhaps more natural than that of $H_\al$. For example, $A_\al$ has some important infinite dimensional modules, which $H_\al$ `cannot see'. 

So it is possible that understanding irreducible $A_\al$-modules is a `more manageable' and more natural task than understanding irreducible $H_\al$-modules. By the way, irreducible $H_\al$-modules $D^\mu$ can be distinguished among all irreducible $A_\al$-modules by the simple condition that all words $\bi=(i_1,\dots,i_n)$ appearing in the formal character of $D^\mu$ have the property that $i_1=0\neq i_2$, see \cite[Proposition~2.4]{LV}. 

One of the interesting ideas, due to Turner \cite{Turner} and others, is `incorporating' representation theory of smaller symmetric groups or Schur algebras into representation theory of $\Si_n$. Curiously, this phenomenon is appearing very naturally in representation theory of $A_\al$ in the form of the so-called imaginary Schur-Weyl duality described below. %, see for example Theorem~\ref{TGGFrag}. 
The results of this section are mainly from \cite{Kcusp} and \cite{KMuth}; a different approach is suggested in \cite{TW}.

\subsection{Convex preorders}
Recall the Cartan matrix $\Car$ and the  simple roots $\al_i$ labeled by  $i\in I$. Let $I'=I\setminus \{0\}$.  
As in \cite{Kac}, we have the affine {\em root system}\, $\Phi$ and the corresponding finite root subsystem $\Phi' =\Phi\cap \Z\,\text{-}\spa(\al_i\mid i\in I')$. Denote by $\Phi'_+$ and $\Phi_+$ the sets of {\em positive}\, roots in $\Phi'$ and $\Phi$, respectively. Then 
$\Phi_+=\Phi_+^\im\sqcup \Phi_+^\re
$, where
$
\Phi_+^\im=\{n\de\mid n\in\Z_{>0}\}
%\qquad\text{and}\qquad \Phi_+^\re=\Phi_{>}^\re\sqcup \Phi_{<}^\re,
$
for the null-root $\de:=\sum_{i\in I}\al_i$. %, and
%$$\Phi_+^\re=\{\be+n\de\mid \be\in  \Phi'_+,\ n\in\Z_{\geq 0}\}\sqcup \{-\be+n\de\mid \be\in  \Phi'_+,\ n\in\Z_{> 0}\}. $$

A {\em convex preorder} on $\Phi_+$ is a preorder $\preceq$ such that the following three conditions hold for all $\be,\ga\in\Phi_+$:
\begin{enumerate}
\item[{\rm (1)}] $\be\preceq\ga \ \text{or}\ \ga\preceq \be$; 
\item[{\rm (2)}] if $\be\preceq \ga$ and $\be+\ga\in\Phi_+$, then $\be\preceq\be+\ga\preceq\ga$;
\item[{\rm (3)}] $\be\preceq\ga$ and $\ga\preceq\be$ if and only if $\be$ and $\ga$ are proportional. 
\end{enumerate}
%(1) $\be\preceq\ga \ \text{or}\ \ga\preceq \be$; 
%(2) if $\be\preceq \ga$ and $\be+\ga\in\Phi_+$, then $\be\preceq\be+\ga\preceq\ga$; 
%(3) $\be\preceq\ga$ and $\ga\preceq\be$ if and only if $\be$ and $\ga$ are proportional. 

We fix a convex preorder $\preceq$ on $\Phi_+$ such that $\al_i\succ n\de\succ\al_0$ for all $i\in I'$; this is always possible. (This additional assumption is for convenient only.) 
Then  
\begin{align*}
\{\be\in \Phi_+^\re\mid \be\succ\de\}&=\{\be+n\de\mid \be\in  \Phi'_+,\ n\in\Z_{\geq 0}\},
\\ 
\{\be\in \Phi_+^\re\mid \be\prec\de\}&=\{-\be+n\de\mid \be\in  \Phi'_+,\ n\in\Z_{> 0}\}. 
\end{align*}
We have  that $\be\preceq\ga$ and $\ga\preceq\be$ happens for $\be\neq\ga$ if and only if both $\be$ and $\ga$ are imaginary. We write $\be\prec\ga$ if $\be\preceq\ga$ but $\ga\not\preceq\be$. 
The following set is {\em totally ordered}\, with respect to $\preceq$:
\begin{equation*}\label{EPsiNew}
\Psi:=\Phi_+^\re\cup\{\de\}.
\end{equation*}

%It is possible to label real roots by non-zero integers so that $\Phi_+^\re=\{\rho_n\mid n\in \Z_{\neq 0}\}$ and 
%\begin{equation}\label{EOrderRoots}
%\Phi^\re_{\succ}=\{\rho_1\succ\rho_2\succ\rho_3\succ\dots\}\ \text{and}\ 
%\Phi^\re_{\prec}=\{\dots\succ\rho_{-3}\succ\rho_{-2}\succ\rho_{-1}\}. 
%\end{equation}

Let $l:=p-1$. 
An {\em $l$-multipartition} of $n$ is a tuple $\umu=(\mu^{(1)},\dots,\mu^{(l)})$ of partitions such that $|\mu^{(1)}|+\dots+|\mu^{(l)}|=n$. The set of all $l$-multipartitions of $n$ is denoted by $\Par_n^l$, and $\Par^l:=\sqcup_{n\geq 0}\Par_n^l$. 
A {\em root partition of $\al\in Q_+$} is a pair $(M,\umu)$, where $M$ is a tuple $(m_\rho)_{\rho\in\Psi}$ of non-negative integers such that $\sum_{\rho\in\Psi}m_\rho\rho=\al$, and $\umu$ is an $l$-multipartition of $m_\de$. 
Clearly all but finitely many integers $m_\rho$ are zero, so we can always choose a finite subset
$$
\rho_1\succ\dots\succ\rho_s\succ\de\succ\rho_{-t}\succ\dots\succ\rho_{-1}
$$
of $\Psi$ such that $m_\rho=0$ for $\rho$ outside of this subset. Then, denoting $m_u:=m_{\rho_u}$, we can 
 write any root partition of $\al$ in the form
\begin{equation}\label{ERP}
(M,\umu)=(\rho_1^{m_1},\dots,\rho_s^{m_s},\umu,\rho_{-t}^{m_{-t}},\dots,\rho_{-1}^{m_{-1}}),
\end{equation}
where all $m_u\in\Z_{\geq 0}$, $\umu\in\Par^l$, and 
$\sum_{u=1}^s m_u\rho_u+|\umu|\de+\sum_{u=1}^{t} m_{-u}\rho_{-u}=\al.$
We write $\Pi(\al)$ for the set of all root partitions of~$\al$. 

Denote by $\Seq$ the set of all finitary  tuples  $M=(m_\rho)_{\rho\in\Psi}\in \Z_{\geq 0}^\Psi$ of non-negative integers, so that a root partition is a pair $(M,\umu)$ with $M\in\Seq$ and $\umu\in\Par^l_{m_\de}$. The left lexicographic order on $\Seq$ is denoted $\leq_l$ and the right lexicographic order on $\Seq$ is denoted $\leq_r$. We will use the following {\em bilexicographic} partial order on $\Seq$: 
$$
M\leq N\qquad\text{if and only if}\qquad M\leq_l N \ \text{and}\ M\geq_r N.
$$
We will use the following partial order on the set $\Pi(\al)$ of root partitions of $\al$:
\begin{equation*}\label{EBilex}
(M,\umu)\leq (N,\unu)\quad \text{if and only if}\quad  M< N\ \text{or $M=N$ and}\ \umu=\unu.
\end{equation*}

\subsection{Cuspidal systems}
Recall from (\ref{EIndFunctor}) and (\ref{EResFunctor}) the functors $\Ind_{\al,\be}$ and $\Res_{\al,\be}$. For $M\in\mod{A_\al}$ and $N\in\mod{A_\be}$, denote $$M\circ N:=\Ind_{\al,\be} M\boxtimes N.$$ 
We also write $M^{\circ n}$ for $M\circ \dots\circ M$ ($n$ times). 

A {\em cuspidal system} (for a fixed convex preorder) is the following data:
\begin{enumerate}
\item[$\bullet$] A {\em cuspidal} irreducible $A_\rho$-module $L_\rho$ assigned to every $\rho\in \Phi_+^\re$, i.e. an irreducible $A_\rho$-module with the following property: if $\be,\ga\in Q_+$ are non-zero elements such that $\rho=\be+\ga$ and $\Res_{\be,\ga}L_\rho\neq 0$, then $\beta$ is a sum of roots less than $\rho$ %(in our convex preorder) 
and $\ga$ is a sum of roots greater than $\rho$. 
\item[$\bullet$] An irreducible {\em imaginary}  $A_{n\de}$-module $L(\umu)$ assigned to every $\umu\in\Par^l_n$, i.e. an irreducible $A_{n\de}$-module with the following property: if $\be,\ga\in Q_+\setminus\Phi_+^\im$ are non-zero elements such that $n\de=\be+\ga$ and $\Res_{\be,\ga}L(\umu)\neq 0$, then $\beta$ is a sum of real roots less than $\de$ and $\ga$ is a sum of real roots greater than $\de$. In addition, it is required that $L(\ula)\not\simeq L(\umu)$ unless $\ula=\umu$. 
\end{enumerate}

%It is proved in \cite{Kcusp} that (for a fixed convex preorder) cuspidal modules exist and are determined uniquely up to an  isomorphism. 

Given a root partition 
$
\pi%=(\rho_1^{m_1},\dots,\rho_s^{m_s},\umu,\rho_{-t}^{m_{-t}},\dots,\rho_{-1}^{m_{-1}})\in\Pi(\al)
$ as in  (\ref{ERP}), %a certain explicit integer $\shift(M,\mu)$ is defined, and then 
set
$ %\begin{equation*}\label{EShift}
\shift(\pi):=\sum_{\rho\in \Phi_+^\re} m_\rho(m_\rho-1)/2\in \Z,
$ %\end{equation*}
and define 
the corresponding {\em (proper) standard module}: 
\begin{equation}\label{EStandIntro}
\Stand(\pi):=q^{\shift(\pi)}L_{\rho_1}^{\circ m_1} \circ\dots\circ L_{\rho_s}^{\circ m_s}\circ  L(\umu)\circ L_{\rho_{-t}}^{\circ m_{-t}}\circ\dots\circ L_{\rho_{-1}}^{m_{-1}}. 
\end{equation}

\begin{Theorem}\label{TCusp} 
For any convex preorder there exists a cuspidal system, unique up to permutation of irreducible imaginary modules.  
%$$\{L_\rho\mid \rho\in \Phi_+^\re\}\cup\{L(\umu)\mid \umu\in\Par^l\}.$$ 
Moreover: 
\begin{enumerate}
\item[{\rm (i)}] For every root partition $\pi$, the standard module  
$
\Stand(\pi)
$ has irreducible head; denote this irreducible module $L(\pi)$. 

\item[{\rm (ii)}] $\{L(\pi)\mid \pi\in \Pi(\al)\}$ is a complete and irredundant system of irreducible $A_\al$-modules up to isomorphism and degree shift.

\item[{\rm (iii)}] For every root partition  $\pi$, we have $L(\pi)^\circledast\cong L(\pi)$.  

\item[{\rm (iv)}] For all root partitions $\pi,\si\in\Pi(\al)$, we have that $[\Stand(\pi):L(\pi)]_q=1$, and $[\Stand(\pi):L(\si)]_q\neq 0$ implies $\si\leq \pi$. 

%\item[{\rm (v)}] For all $(M,\umu),(N,\unu)\in\Pi(\al)$, we have that $\Res_{|M|}L(M,\umu)\simeq L_{M,\umu}$ and $\Res_{|N|}L(M,\umu)\neq 0$ implies $N\leq M$.  

\item[{\rm (v)}] The induced module $L_\rho^{\circ n}$ is irreducible for all $\rho\in\Phi^\re_+$ and $n\in\Z_{>0}$.
\end{enumerate}
\end{Theorem}

\subsection{Minuscule representations and imaginary tensor spaces}
Theorem~\ref{TCusp} gives a `rough classification' of irreducible $A_\al$-modules. The main problem is that we did not give a canonical definition of individual irreducible imaginary modules $L(\umu)$. So far, we just know that the amount of such modules for $A_{n\de}$ is equal to the number of $l$-multipartitions of $n$, and 
we have labeled them by such multipartitions in an arbitrary way. 

To address this problem, we begin with an explicit description of the {\em minuscule} representations---the irreducible imaginary $A_\de$-modules. These correspond to $l$-multipartitions of~$1$. There are of course exactly $l$ such multipartitions, namely $\umu(1),\dots,\umu(l)$, where 
$$
\umu(i):=(\emptyset,\dots,\emptyset,(1),\emptyset,\dots,\emptyset)%\qquad(1\leq i\leq l),
$$
with the partition $(1)$ in the $i$th position. 

Let $i\in I'=\{1,2,\dots,p-1\}$ (we identify the set $I$ of residues modulo $p$ with integers $0,1,\dots,p-1$). Consider the hook partition $\chi^i=(i,1^{p-i})$ for all $i\in I'$. For example, if $p=5$, here are the corresponding Young diagram with the residues of the boxes written in them.   
$$
\begin{picture}(-9,44)
\put(-162,0){$\chi^1=$}
\put(-137,-34){$\Tableau{{0},{4},{3},{2},{1}}$}
\put(-100,0){$\chi^2=$}
\put(-72,-25){$\Tableau{{0,1},{4},{3},{2}}$} 
\put(-25,0){$\chi^3=$}
\put(1,-18){$\Tableau{{0,1,2},{4},{3}}$} 
\put(60,0){$\chi^4=$}
\put(90,-11){$\Tableau{{0,1,2,3},{4}}$}
\end{picture}
$$

\vspace{10mm}
Note that the partitions $\chi^i$ for $i\in I'$ are homogenous in the sense of \S\ref{SSHomog}. In particular, we have the corresponding homogeneous irreducible $H_\de$-modules $D^{\chi^i}$ defined explicitly in Theorem~\ref{Thomog}. Define the $A_\de$-modules 
$$
L(\umu(i)):=L_{\de,i}:= \infl D^{\chi^i}\qquad(i\in I'). 
$$
For example, $L_{\de,1}$ and $L_{\de,p-1}$ are $1$-dimensional with characters
$$
\CH L_{\de,1}=(0,p-1,p-2,\dots,1),\quad \CH L_{\de,p-1}=(0,1,2,\dots, p-1),
$$ 
while for $p>3$, the module $L_{\de,p-2}$ is $(p-2)$-dimensional with character
$$
\CH L_{\de,p-2}=\textstyle\sum_{r=0}^{p-3}(0,1,\dots,r, p-1,r+1,\dots,p-2).
$$
Define the {\em imaginary tensor space of color $i\in I'$}  to be the $A_{n\de}$-module
$$
M_{n,i}:=L_{\de,i}^{\circ n}. 
$$

Fix for now $i\in I'$ and suppress $i$ from the indices, so that we have the imaginary tensor space $M_n=M_{n,i}$. 
The $A_{n\de}$-module structure on $M_n$ yields an algebra homomorphism $A_{n\de}\to \End_F(M_n)$. Define the {\em imaginary Schur algebra} $\ImS_n
$ as the image of $A_{n\de}$ under this homomorphism, i.e.  $\ImS_n=A_{n\de}/\Ann_{A_{n\de}}(M_n).$ 
Modules over $A_{n\de}$ which factor through to $\ImS_n$ are called {\em imaginary modules} (of color $i$). It turns out that this notion agrees  with the notion of an irreducible imaginary module in the sense of cuspidal systems.  
%Thus the {\em category of imaginary $A_{n\de}$-modules} is the same as the category of $\ImS_n$-modules. It is clear that $M_n$ and its composition factors are imaginary modules. Conversely, any irreducible $\ImS_n$-module appears as  a composition factor of $M_n$. 

\begin{Theorem} \label{TISWD}%{\rm \cite{}}%{\bf ()}
Let $n\in\Z_{>0}$. Then:
\begin{enumerate}
\item[{\rm (i)}] $M_n$ is a projective $\ImS_n$-module. 
\item[{\rm (ii)}] The endomorphism algebra $\End_{A_{n\de}}(\Mde_n)^{\op}=\End_{\ImS_n}(\Mde_n)^{\op}$ is isomorphic to the group algebra $F\Si_n$ of the symmetric group $\Si_n$ (concentrated in degree zero). Thus $M_n$ can be considered as a right $F\Si_n$-module. 
\item[{\rm (iii)}] $\End_{F\Si_n}(M_n)=\ImS_n$.

\end{enumerate}
\end{Theorem}

In view of the theorem, we have an exact functor 
\begin{equation*}\label{EGaFun}
\ga_{n}: \mod{\ImS_n}\to\mod{F\Si_n},\quad V\mapsto \Hom_{\ImS_n}(M_n,V).
\end{equation*}
Unless $p>n$ or $p=0$, the $\ImS_n$-module $M_n$ is not a projective generator, and $\ga_n$ is not an equivalence of categories. To fix this problem, we need to upgrade from the  {\em imaginary Schur-Weyl duality} of Theorem~\ref{TISWD} to an {\em imaginary Howe duality}.

\subsection{Imaginary How and Ringel dualities}
Let 
$
{\mathtt x}_{n}:=\sum_{g\in \Si_n}g.
$
Define the  {\em imaginary exterior} and {\em  %symmetric, 
imaginary divided} powers respectively as follows:
\begin{align*}
%\Sde_n&:=\Mde_n/\spa\{mg-\sgn(g)m\mid g\in\Si_n,\ m\in \Mde_n\},\\
\Lade_n:=\Mde_n {\mathtt x}_n,
\quad
\Zde_n:=\{m\in \Mde_n\mid mg-\sgn(g)m=0\ \text{for all $g\in\Si_n$}\}.
\end{align*}

For $h\in\Z_{>0}$, denote by $X(h,n)$ the set of all compositions of $n$ with $h$ parts:
$$
X(h,n):=\{(n_1,\dots,n_h)\in\Z_{\geq 0}^h\mid n_1+\dots+n_h=n\}.
$$
The corresponding set of partitions is
$$
X_+(h,n):=\{(n_1,\dots,n_h)\in X(h,n)\mid n_1\geq\dots\geq n_h\}.
$$
For a composition $\nu=(n_1,\dots,n_h)\in X(h,n)$, we define the functor of {\em imaginary induction}: %  and {\em imaginary restriction}:
$$
\HCI_\nu^n:=\Ind_{n_1\de,\dots,n_h\de}: \mod{A_{n_1\de,\dots,n_h\de}}\to \mod{A_{n\de}}.
$$
%and$$\HCR_\nu^n:=\Res_{n_1\de,\dots,n_h\de}: \mod{A_{n\de}}\to \mod{A_{n_1\de,\dots,n_h\de}}.$$
%These functors `respect' the categories of imaginary representations. For example, g
Given imaginary $A_{n_b\de}$-modules $V_b$ for $b=1,\dots,h$, the module $\HCI_\nu^n(V_1\boxtimes\dots\boxtimes V_h)$ is also imaginary. Define
\begin{align*}
%\Sde^\nu&:=\HCI_\nu^n(\Sde_{n_1}\boxtimes\dots\boxtimes \Sde_{n_h}),\\
\Zde^\nu:=\HCI_\nu^n(\Zde_{n_1}\boxtimes\dots\boxtimes \Zde_{n_h}),\quad 
\Lade^\nu:=\HCI_\nu^n(\Lade_{n_1}\boxtimes\dots\boxtimes \Lade_{n_h}).
\end{align*} 

Now, let $S_{h,n}$ be the classical Schur algebra, whose representations are the same as the degree $n$ polynomial representations of the general linear group $GL_h(F)$ \cite{Green}. It is a finite dimensional quasi-hereditary algebra with irreducible, standard, costandard, and indecomposable tilting modules
$$
L_h(\la),\ \De_h(\la),\ \nabla_h(\la),\ T_h(\la)\qquad(\la\in X_+(h,n)). 
$$

\begin{Theorem} %\label{}%{\rm \cite{}}%{\bf ()}
We have: 
\begin{enumerate}
\item[{\rm (i)}] For each $\nu\in X(h,n)$ the $\ImS_n$-module $\Zdot^\nu$ is  projective. Moreover, for any $h\geq n$, the module $Z:=\bigoplus_{\nu\in X(h,n)}\Zdot^\nu$ is a projective generator for $\ImS_n$.  
\item[{\rm (ii)}] The endomorphism algebra $\End_{\ImS_n}(Z)^{\op}$ is isomorphic to the classical Schur algebra $S_{h,n}$ concentrated in degree zero. Thus $Z$ can be considered as a right $S_{h,n}$-module. 
\item[{\rm (iii)}] $\End_{S_{h,n}}(Z)=\ImS_n$. 
\end{enumerate}
\end{Theorem}

This theorem allows us to use Morita theory and define quasi-inverse  equivalences of categories: 
\begin{align}
\al_{h,n}&: \mod{\ImS_n}\to\mod{S_{h,n}},\quad V\mapsto \Hom_{\ImS_n}(Z,V)
\label{EAlpha}
\\
\be_{h,n}&:\mod{S_{h,n}}\to \mod{\ImS_n},\quad W\mapsto Z\otimes_{S_{h,n}}W.
\label{EBeta}
\end{align}

Let $\mu \in \Par_n$ and $h\geq n$. We can also consider $\mu$ as an element of $X_+(h,n)$. Define the $\ImS_n$-modules: %(hence, by inflation, also graded $A_{n\de}$-modules):
\begin{align*} %\begin{eqnarray*}
\Lde(\mu) &:= \be_{h,n}(L_h(\mu)),
\ 
\Dede(\mu) := \be_{h,n}(\Delta_h(\mu)),
\\ 
\nade(\mu) &:= \be_{h,n}(\nabla_h(\mu)),
\ 
T(\mu) := \be_{h,n}(T_h(\mu)).
\end{align*} % \end{eqnarray*}
These definitions turn out to be independent of the choice of $h \geq n$. %, see \cite[Lemma~6.1.3]{KMuth}. 
An easy consequence of the theorem above is that the imaginary Schur algebra $\ImS_n$ is a finite dimensional quasi-hereditary algebra with irreducible, standard, costandard, and indecomposable tilting modules
$
L(\mu),\ \De(\mu),\ \nabla(\mu),\ T(\mu) 
$ 
labeled by $\mu\in X_+(h,n)$. In particular, inflating the irreducible modules $L(\mu)$ from $\ImS_n$ to $A_{n\de}$, we get:

\begin{Theorem} \label{TIrrOneColor}%{\rm \cite{}}%{\bf ()}
The irreducible imaginary $A_{n\de}$-modules of color $i$ are exactly the modules $\{L(\la)\mid\la\in \Par_n\}$ (up to isomorphism).  
\end{Theorem}

Moreover:

\iffalse{
Denoting by $f_{h,n}$ the usual Schur functor, as for example in \cite[\S6]{Green}, by definitions we then have a commutative triangle up to isomorphism of functors:
\begin{equation*}\label{EFunTriangle}
\begin{picture}(100,30)
\put(0,-30){
\begin{tikzpicture}%[scale=1, line join=bevel]
\node at (0,0.5) {$\mod{S_{h,n}}$};
\node at (1.5,-1) {$\mod{\ImS_n}$}; 
\node at (-1.5,-1) {$\mod{F\Si_n}$};

\draw [->] (-0.5,0.2) -- (-1.2,-0.8);
%\draw [-] (1.6,1.2) -- (1.6,0.3);
\draw [<-] (-0.5,-1) -- (0.6,-1);
\draw [->] (0.5,0.2) -- (1.3,-0.7);
\draw [<-] (0.3,0.2) -- (1.1,-0.7);

\node at (-1.3,-0.2) {${\scriptstyle f_{h,n}}$};
\node at (0.25,-0.25) {${\scriptstyle \al_{h,n}}$};
\node at (1.35,-0.25) {${\scriptstyle \be_{h,n}}$};
\node at (0,-0.8) {${\scriptstyle \ga_{n}}$};
\end{tikzpicture}
}
\end{picture}
\end{equation*}

\vspace{1cm}

}\fi

%Since the Schur functor sends tensor products to the induction for symmetric groups, cf. \cite[Theorem 4.13]{BKlr}, we deduce from Theorem 5 and (\ref{EFunTriangle}):

\iffalse{

\vspace{2mm}
\noindent
{\bf Corollary.}
{\em
Let $\nu=(n_1,\dots,n_a)\in X(a,n)$. We have a functorial isomorphism  
$$
\ga_n\big(\HCI_\nu^n (V_1\boxtimes\dots\boxtimes  V_a)\big)\cong \ind_{F\Si_\nu}^{F\Si_n}
\big(\ga_{n_1}(V_1)\boxtimes\dots\boxtimes \ga_{n_a}(V_a)\big)
$$
for $V_1\in\mod{\ImS_{n_1}},\dots,V_a\in\mod{\ImS_{n_a}}$.
}

\vspace{2mm}
}\fi

\begin{Theorem} %\label{}%{\rm \cite{}}%{\bf ()}
We have:
\begin{enumerate}
\item[{\rm (i)}] Let $h\geq n$. The $\ImS_n$-module $\bigoplus_{\nu\in X(h,n)} \Lade^\nu$ is a full tilting module. 
\item[{\rm (ii)}] We have isomorphisms of endomorphism algebras 
$$\textstyle\End_{\ImS_n}\big(\bigoplus_{\nu\in X(h,n)} \Lade^\nu \big)^{\op}\cong S_{h,n}\quad\text{and}\quad   \End_{S_{h,n}}\big(\bigoplus_{\nu\in X(h,n)} \Lade^\nu \big)\cong \ImS_{n}.$$
\end{enumerate}
\end{Theorem}

The additional nice property of the constructed Morita equivalence is that  imaginary induction commutes with tensor products:

\begin{Theorem} %\label{}%{\rm \cite{}}%{\bf ()}
Let $h\geq n$ and $\nu=(n_1,\dots,n_a)\in X(a,n)$. 
The following functors are isomorphic:
\begin{align*}
\HCI_\nu^n (\be_{h,n_1}-\boxtimes\dots\boxtimes \be_{h,n_a} -) & :\mod{S_{h,n_1}}\times\dots\times \mod{S_{h,n_a}}\to \mod{\ImS_n},
\\
\be_{h,n}(-\otimes\dots\otimes  -) & :\mod{S_{h,n_1}}\times\dots\times \mod{S_{h,n_a}}\to \mod{\ImS_n}.
\end{align*}

\end{Theorem}

\subsection{Gelfand-Graev character fragment and imaginary Jacobi-Trudi formula}

We can say quite a bit about the characters of irreducible imaginary modules. An important role in the paper is played by an analogue of the {\em Gelfand-Graev representation}, cf. e.g. \cite{BDK}. %We give its definition and describe some of the key properties. 
Let $$\bi=(i_1,\dots, i_p)$$ be any word appearing in the (explicitly known) formal character of $L_\de$. Define the corresponding {\em Gelfand-Graev  words} 
$
\GG^{(t)}_\bi:=i_1^t i_2^t\dots i_p^t 
$ 
for all $t\in\Z_{>0}$.
For any composition 
$\mu=(\mu_1,\dots,\mu_n)\in X(n, n)$
define 
$$ %\begin{equation}\label{EGGW}
\GG^\mu_\bi:=\GG_\bi^{(\mu_1)}\dots\GG_\bi^{(\mu_n)} 
$$ %\end{equation}
 and 
$ %\begin{equation}\label{EmlaIntro}
c_\bi(\mu):=([\mu_1]_q^!\dots [\mu_n]_q^!)^p\in\Laurent.  
$ %\end{equation}
If $V\in\mod{A_{n\de}}$,  it is known that   
$$\DIM V_{\GG_\bi^\mu}=c_\bi(\mu)m_{\bi,\mu}(V)$$ 
for some $m_{\bi,\mu}(V)\in\Laurent$. We are going to describe $m_{\bi,\mu}(V)$ for many important imaginary modules.

There are explicitly defined {\em Gelfand-Graev idempotents}\, $\ga_{n,\bi}\in A_{n\de}$. 
The {\em Gelfand-Graev module} is the projective module 
$
\Ga_{n,\bi}:=q^{-pn(n-1)/2} A_{n\de}\ga_{n,\bi}.
$

\begin{Theorem} %\label{}%{\rm \cite{}}%{\bf ()}
%Let $\bi$ be any word appearing in the formal character of $L_\de$. 
For any $V\in\mod{A_{n\de}}$ and $\mu\in X(n,n)$, we have 
$$
m_{\bi,\mu}(V)=\DIM\Hom_{A_{n\de}}(\Ga_{\mu_1,\bi}\circ \dots\circ \Ga_{\mu_n,\bi},V).
$$
\end{Theorem}

Since we have equivalences of categories (\ref{EAlpha}) and (\ref{EBeta}), every finite dimensional graded $\ImS_n$-module $V$ can be written as $V=\be_{n,n}(W)$ (up to degree shift). Then we can describe the {\em Gelfand-Graev fragment} of $\CH V$ as follows:

\begin{Theorem} \label{TGGFrag}%{\rm \cite{}}%{\bf ()}
Let $\bi$ be any word appearing in the formal character of $L_\de$, %$\la\in X_+(n,n)$, 
$\mu\in X(n, n)$, $W\in\mod{S_{h,n}}$, and $V=\beta_{n,n}(W)\in\mod{\ImS_n}$. Denote by $W_\mu$ the usual weight space of $W$. Then 
$$\DIM V_{\GG_\bi^\mu}=c_\bi(\mu)\dim W_\mu.$$
\end{Theorem}

Note that the Gelfand-Graev fragment is described in terms of the formal characters of a `smaller rank' Schur algebra.

The formal characters of the modules $\De(\mu)$ are important; for example in the case $p>n$ we have $\De(\mu)=L(\mu)$. An {\em imaginary Jacobi-Trudi formula}\, allows us to compute the formal characters of the modules $\De(\la)$ explicitly. 

First of all, the characters of the modules $\De(1^n)=\La_n=L(1^n)$ are well-understood: let $\bi=(i_1,\dots,i_p)$ be a word appearing in $L_\de$. Then $\bi^n$ is a word of $\De(1^n)$, and $\De(1^n)$ is the homogeneous irreducible module associated to the connected component of $\bi^n$ in the word graph, see \cite{KRhomog}. 

Let $\mu=(\mu_1,\dots,\mu_a)\in \Par_n$. Denote by $\circ$ the quantum shuffle product, see e.g. \cite[\S2]{KL1}. Then $\CH\De(1^k)\,\circ\,\CH\De(1^l)=\CH\De(1^l)\,\circ\, \CH\De(1^k)$ for all $k,l\in\Z_{>0}$. 
So we can use the quantum shuffle product to make sense of the  following determinant as an element of $\A \words_{n\de}$: 
$$
\JTD(\mu):=\det\big(\CH\De(1^{\mu_r-r+s})\big)_{1\leq r,s\leq a}.
$$
where $\CH\De(1^{0})$ is interpreted as (multiplicative) identity, and $\CH\De(1^{m})$ is interpreted as (multiplicative) zero if $m<0$. 
For example, for $\mu=(3,1,1)$:
\begin{align*}
\JTD((3,1,1))
=&\det 
\left(
\begin{matrix}
\CH\De(1^{3}) & \CH\De(1^{4})& \CH\De(1^{5})  \\
1 & \CH\De(1)& \CH\De(1^{2})\\
0 & 1 & \CH\De(1)
\end{matrix}
\right)
\\
=&\CH\De(1^{3})\circ\CH\De(1)\circ \CH\De(1)+\CH\De(1^{5})
\\&-\CH\De(1^{4})\circ\CH\De(1)-\CH\De(1^{3})\circ \CH\De(1^{2}).
\end{align*}

\begin{Theorem} \label{TDet} %{\rm \cite{}}% 
Let $\mu^\tr$ be the partition transpose to $\mu$. 
Then $
\CH\Dede(\mu)=\JTD(\mu^\tr)$. 
\end{Theorem}

For example, let $p=2$. Then $I'=\{1\}$. The $i=1$. In this case the character of $L_\de$ is $(0,1)$, and the character of $\De(1^n)$ is $(0,1,0,1,\dots,0,1)$. So 
$$
D((1,1))=\det 
\left(
\begin{matrix}
\bi_{\de} & \bi_{2\de}  \\
1 & \bi_{\de}
\end{matrix}
\right)
=\bi_{\de}\circ\bi_{\de}-\bi_{2\de}=(0101)+(q+q^{-1})^2 (0011). 
$$

\subsection{Classification of imaginary irreducible modules}
In Theorem~\ref{TIrrOneColor}, we have classified the irreducible imaginary representations of $A_{n\de}$ of color $i$. Since we now want to distinguish between the imaginary representations of different colors, we will use the notation $L_i(\mu)$ for these irreducible imaginary representations of color $i$ corresponding to a partition $\mu$.

\begin{Theorem} %\label{}%{\rm \cite{}}%{\bf ()}
For an $l$-multipartition $\ula=(\la^{(1)},\dots,\la^{(l)})$ of $n$, define $$L(\ula):=L_1(\la^{(1)})\circ\dots\circ L_l(\la^{(l)}).$$ 
Then $\{L(\ula)\mid\ula\in\Par^l_n\}$ is a complete and irredundant system of imaginary irreducible  $A_{n\de}$-modules. 
\end{Theorem}

\end{document}